\documentclass[letterpaper]{article}
\usepackage{mathrsfs}
\usepackage{amssymb}
\usepackage{amsfonts}
\usepackage{graphicx}
\usepackage{subfig}
\usepackage{float}
\usepackage[titletoc]{appendix}
\usepackage{array}
\usepackage{amsmath}
\usepackage{epsf}
\usepackage{indentfirst}
\usepackage{amsthm}
\usepackage{color}

\usepackage{hyperref}
\graphicspath{{../}}

\newcommand{\bd}{\boldsymbol}

\newcommand{\pl}[2]{\frac{\partial#1}{\partial#2}}
\newcommand{\p}{\partial}
\newcommand{\f}[2]{\frac{#1}{#2}}

\newcommand{\nn}{\nonumber}

\newcommand{\dt}{\delta}

\newcommand{\ld}{\lambda}
\newcommand{\Ld}{\Lambda}

\newcommand{\Tht}{\Theta}

\newcommand{\vep}{\varepsilon}

\newcommand{\ov}{\overline}
\newcommand{\be}{\begin{equation}}
\newcommand{\ee}{\end{equation}}
\newcommand{\ba}{\begin{array}}
\newcommand{\ea}{\end{array}}
\newcommand{\bea}{\begin{eqnarray}}
\newcommand{\eea}{\end{eqnarray}}
\newcommand{\beas}{\begin{eqnarray*}}
\newcommand{\eeas}{\end{eqnarray*}}
\newcommand{\bl}{\begin{law}}
\newcommand{\el}{\end{law}}
\newcommand{\bthm}{\begin{thm}}
\newcommand{\ethm}{\end{thm}}

\newcommand{\mbb}{\mathbb}
\newcommand{\mcl}{\mathcal}

\newcommand{\e}{\mathrm{e}}
\newcommand{\ri}{\mathrm{i}}
\newcommand{\rd}{\mathrm{d}}
\newcommand{\R}{\mathbb{R}}

\newcommand{\grdx}{\nabla_{\bd{x}}}
\newcommand{\grdp}{\nabla_{\bd{p}}}
\newcommand{\lapx}{\Delta_{\bd{x}}}

\begin{document}

\title{A multi-band semiclassical model for surface hopping quantum dynamics\thanks{This work was partially supported by NSF grants DMS-1114546 and DMS-1107291: NSF Research Network in Mathematical Sciences KI-Net: Kinetic description of emerging challenges in multiscale problems of natural sciences.}}
\date{}
\author{
Lihui Chai\thanks{Department of Mathematics, University of California, Santa Barbara, CA 93106, USA (chailh@math.ucsb.edu)},
Shi Jin\thanks{Department of Mathematics, Institute of Natural Sciences and
MOE-LSE, Shanghai Jiao Tong University, Shanghai
200240, China; and Department of Mathematics, University of Wisconsin-Madison, Madison, WI 53706, USA
(jin@math.wisc.edu)},
Qin Li\thanks{Department of Computing and Mathematical Sciences,
California Institute of Technology, Pasadena, CA 91125, USA
(qinli@caltech.edu)}, and
Omar Morandi\thanks{Institut de Physique et Chimie des Mat\'eriaux de Strasbourg, CNRS and University of Strasbourg,  
23, rue du Loess, F-67034 Strasbourg, France (omar.morandi@ipcms.unistra.fr)}
}

\maketitle

\numberwithin{equation}{section} \numberwithin{figure}{section}

\newtheorem{law}{Lemma}[section]
\newtheorem{thm}{Theorem}[section]
\theoremstyle{definition}\newtheorem{remark}{Remark}[section]
\theoremstyle{definition}\newtheorem{lemma}{Lemma}[section]
\theoremstyle{definition}\newtheorem{corollary}{Corollary}[section]
\theoremstyle{definition}\newtheorem{assumption}{Assumption}

%\tableofcontents
\begin{abstract}
In the paper we derive a semiclassical model for surface hopping allowing
quantum dynamical non-adiabatic transition between different potential
energy surfaces in which cases the classical Born-Oppenheimer approximation
breaks down. The model is derived using the Wigner transform and
Weyl quantization, and the central idea is to evolve the {\it entire} Wigner
matrix rather than just the diagonal entries as was done previously
in the adiabatic case.
The off-diagonal entries of the Wigner matrix
suitably describe the non-adiabatic transition, such as the Berry connection, for avoided crossings.
We study the numerical approximation issues of the model,
and then conduct numerical experiments to validate the model.
\end{abstract}

\section{Introduction}
In this paper we derive a semiclassical model based on the quantum phase-space description of the particle dynamics. We consider
 the nucleonic Schr\"odinger system:
 \begin{align}\label{eqn_Schrodinger}
  \ri\vep\pl{\psi^\vep}{t}(t,\bd{x}) &= \hat{H}\psi^\varepsilon(t,x),\qquad\left(t,\bd{x}\right)\in\left(\mbb{R}^+,\mbb{R}^d\right)\\
  \psi^\varepsilon(0,\bd{x})&=\psi^\varepsilon_0(\bd{x})
 \end{align}
 with the self-adjoint Hamiltonian operator defined by:
\begin{align}\label{eqn_Hamiltonian}
 \hat{H} = -\frac{\vep^2}{2}\Delta_{\bd{x}}+\tilde{V}(\bd{x}).
\end{align}
Here, $\psi$ is a vector and $\tilde{V}(x)$ is a Hermitian matrix.
$\varepsilon = \sqrt{\frac{m}{M}}$ is the mass ratio between electron and
nuclei. This system arises from the Born-Oppenheimer approximation \cite{BO}
of the $N-$body Schr\"odinger equation in which the nucleonic Schr\"odinger system
\eqref{eqn_Schrodinger} is solved along the electronic potential surfaces. We will focus on the two-energy system although
our study can be extended to systems with more energy levels in a straightforward way. In the two-energy level case, the potential matrix reads as:
\begin{align}\label{eqn_potential}
\begin{split}
 \tilde{V}(\bd{x})&=\frac{1}{2}\text{tr}\tilde{V}(\bd{x})+V=U(\bd{x})+
 \left(\begin{array}{cc}u(\bd{x})
 & v(\bd{x})\\v^\dagger(\bd{x}) & -u(\bd{x})\end{array}\right).
\end{split}
\end{align}

For future references, we consider the unitary matrix 
 $\Theta$
\begin{equation}\label{eqn_theta}
 \Theta^\dagger(\bd{x})=\left[\chi_+,\chi_-\right],
\end{equation}
that diagonalizes the potential operator $V$. We have
\begin{align}\label{eqn_diagonalize}
 V=\Theta^\dagger\Ld_V\Theta,
\end{align}
where
\begin{equation}\label{eqn_lambda}
 \Ld_V(\bd{x})=\mathrm{diag}\big(E(\bd{x}),\,-E(\bd{x})\big)
 =\mathrm{diag}\left(\sqrt{|u(\bd{x})|^2+|v(\bd{x})|^2},-\sqrt{|u(\bd{x})|^2+|v(\bd{x})|^2}\right),
\end{equation}
and
\begin{align}\label{eqn_Tht}
\begin{split}
\Tht&=\frac{1}{ \sqrt{2\,\left(1+\f{u(\bd{x})}{E(\bd{x})}\right)}}
    \left(\begin{array}{rr}
    \left(1+\f{u(\bd{x})}{E(\bd{x})}\right)\f{v^\dagger(\bd{x})}{|v(\bd{x})|} &
    \f{|v(\bd{x})|}{E(\bd{x})}\\
    -\f{v^\dagger(\bd{x})}{E(\bd{x})} &
    1+\f{u(\bd{x})}{E(\bd{x})}
    \end{array}\right).
\end{split}
\end{align}
Obviously $\chi_\pm$ are the eigenvectors of $V$ corresponding to the eigenvalues $\pm E$ with  $E(\bd{x})=\sqrt{|u(\bd{x})|^2+|v(\bd{x})|^2}$.  Hereafter, we call the two eigenvalues the {\it energy bands}, and $\Delta{E}=2E$, the {\it energy gap}.

For $\Delta E=0$, the matrix $\Theta$ becomes singular (conical crossing). In this paper, we are interested in the cases where the energy gap is strictly positive and asymptotically small. In particular, we focus on the so called  {\it avoided crossing} scaling where the minimum of the energy gap is of the order $\sqrt{\varepsilon}$.

We consider a few types of prototype potentials and analyze their influence on the non-adiabatic transitions process:
\begin{description}
\item[1D case:]
 \begin{equation}\label{eqn_u_v_delta_1d}
  u(x)=x,\quad v(x,\dt)\equiv \delta,\quad U(\bd{x}) = 0.
 \end{equation}
 The eigenvalues are $\Ld_V^{\pm}=\pm E=\pm\sqrt{x^2+\dt^2}$ and the avoided crossing point is $x=0$.
\item[2D cases:]
 First example
 \begin{equation}\label{eqn_u_v_delta_real}
  u(\bd{x},\dt)=x,\quad v(\bd{x},\dt)=\sqrt{y^2+\dt^2},\quad U(\bd{x})=0.
 \end{equation}
 Second example
 \begin{equation}\label{eqn_u_v_delta}
  u(\bd{x},\dt)=x,\quad v(\bd{x},\dt)=y+\ri\dt,\quad U(\bd{x})=0.
 \end{equation}
 Here, we denoted $\bd{x}=(x,y)$. In the 2D cases, the eigenvalues are given by  $\Ld_V^{\pm}=\pm E=\pm\sqrt{x^2+y^2+\delta^2}$, and the avoided crossing  point is $(x,y)=(0,0)$.
\end{description}

We are interested in deriving a semiclassical approximation to the
Schr\"odinger system (\ref{eqn_Schrodinger}) with avoided-crossings. One of the advantages of our method is that the computational cost is significantly reduced than solving directly the original Schr\"odinger system \eqref{eqn_Schrodinger}.

One of the difficulties in the derivation of the semiclassical expansion for a system with two or more energy levels, is the non-commutativity of the matrix $\tilde{V}$ with the Laplacian operator. In the case where $\Delta{E}$ is of order 1, the equation of motion can be well approximated by a fully diagonalized system -- one classical Liouville equation for each energy level~\cite{kato_adiabatic_1950,joye_exponential_1991, hagedorn_time_1980,hagedorn_high_1986,hagedorn_high_1988,hagedorn_time-dependent_2001, avron_adiabatic_1998,panati_space-adiabatic_2002,spohn_adiabatic_2001,nenciu_semiclassical_2004}. This is the standard Born-Oppenheimer approximation. See reviews \cite{hagedorn_mathematical_2007,teufel_adiabatic_book}. However,  when the eigenvalues are of $O(\sqrt{\varepsilon})$ away from each other, the classical Born-Oppenheimer approximation breaks down and the diagonalized system is no longer a good approximation of the full coupled system. In such case, around the crossing points, the particles 
could move from one band to another (the {\it non-adiabatic} phenomenon).

The study of the mathematical properties and the study of the physical systems where the energy band structure shows some crossing points dates back to Wigner and von Neumann~\cite{von_neuman_uber_1929}. It can be shown that the crossing set is of measure zero, while the influence is of order $1$, and it is this crossing phenomenon that is responsible of some chemical reactions~\cite{worth_beyond_2004,yarkony_nonadiabatic_2012}. Due to its physical significance, this topic has been studied extensively in computational chemistry community. The first result on the transition rate is due to Landau and Zener~\cite{zener_non-adiabatic_1932}, who gave a rough estimate on the transition probability. Afterwards, there is very rich literature investigating the different aspects of the problem, including theoretical studies and algorithm development. We here mention the two most well-known algorithms, both by Tully {\it etc.}: the surface hopping method based on applying Landau-Zener formula~\cite{tully_trajectory_1971}, and the fewest switches method~\cite{tully_molecular_1990}, a Markov-Chain Monte-Carlo type method. Some criticisms have been also raised on the Landau-Zener formula, and we mention~\cite{bates_collisions_1960,herman_generalization_1982}.

On the mathematical side, in~\cite{hagedorn_proof_1991} Hagedorn firstly rigorously reexamined Zener's idea. This was followed by a series of works~\cite{joye_exponentially_1991,joye_exponential_1991,jaksic_landauzener_1993,hagedorn_molecular_1999,hagedorn_landau-zener_1998,rousse_landau-zener_2004,bourquin_non-adiabatic_2011}, in which they also show that the jumping behavior could heavily depends on the types of crossings (see classification of crossings in~\cite{hagedorn_classification_1998}). The study of the non-adiabatic transition on the phase space was done in~\cite{fermanian_kammerer_wigner_2003}.

The surface hopping algorithms that use the Landau-Zener formula for evaluating the non-adiabatic transitions for conical crossings have seen recent mathematical interests~\cite{lasser_construction_2007,jin_eulerian_2011,jin_hybrid_2011,fermanian-kammerer_single_2011}. The main advantage of these surface hopping methods, compared to computing the original Schr\"odinger equation~\eqref{eqn_Schrodinger}, is that they do not need to numerically resolve the $O(\varepsilon)$ wavelength. However, these methods cannot account for phase information at the crossing points, and thus ignore important physical phenomenon~\cite{Niu,Morandi2011}. The main result of this paper is to present a semiclassical model that includes the particle phase correction at the crossing points.

Our method is based on the Wigner transform~\cite{lions_sur_1993} and the Weyl quantization~\cite{weyl_quantenmechanik_1927,hormander_weyl_1979} procedure. In the adiabatic case, with $\varepsilon\to 0$ (classical limit), the Wigner transformation leads to a set of decoupled Liouville equations, each for one energy band \cite{lions_sur_1993, gerard_homogenization_1997}. In this case, only the diagonal entries of the Wigner matrix that correspond to the projection onto the two eigenspaces of the underlying Hamiltonian are relevant. However, in the presence of a band crossing, one cannot ignore the off-diagonal terms. For this reason, our main idea in the paper is to find the semiclassical approximation for the {\it entire Wigner matrix}. This approach is similar to the derivation of the transport equation for graphene \cite{morandi_wigner_2011} and in semiconductors systems \cite{Morandi2009PRB}. Our model is a coupled Liouville system for all entries of the Wigner matrix, where the off-diagonal terms 
prescribes the quantum transition between bands,  and the two-bands correlations due to the Berry connections. We also discuss numerical approximation of this model utilizing a multi-physics domain decomposition idea proposed in~\cite{chai_semi-classical_2013}: away from the crossing points we solve the standard adiabatic Liouville equations, while in the crossing zones the new semiclassical system is solved, and the two systems are connected by interface conditions.

In section~\ref{sec_model} we present the derivation of the new semiclassical model. We also produce a primitive analysis of the behavior of the solutions to the system. In section~\ref{sec_NumDDM} we describe a coupling method that combines the new semiclassical model near the crossing points with the adiabatic system elsewhere in order to further reduce the computational cost. Numerical examples are shown afterwards.

\section{The semiclassical formulation}\label{sec_model}
In the following we describe the basics of the Wigner transform and the Weyl quantization. In subsections \ref{sec_Adiabatic} and~\ref{sec_Non-adiabatic} we derive the mathematical model for the adiabatic and non-adiabatic cases respectively.

\subsection{The Wigner transformation and Weyl quantization}\label{sec_Wigner}

The Wigner function is defined by
\begin{equation}\label{def_WignerTransform}
 F^\varepsilon(\bd{x},\bd{p})=\f{1}{(2\pi)^d}\int\,\rho^\varepsilon\left(\bd{x}-\frac{\vep\bd{y}}{2},\bd{x}+\frac{\vep\bd{y}}{2}\right)\,\e^{{\ri}\bd{p}\cdot\bd{y}}\,\rd\bd{y},
\end{equation}
where $\rho^\varepsilon(\bd{x},\bd{x}')=\psi^\vep(\bd{x})\otimes\ov{\psi^\vep}(\bd{x}')$ is the density matrix, $\psi^\vep$ is defined in \eqref{eqn_Schrodinger}, $\ov{\psi^\varepsilon}$ is the complex conjugate of $\psi^\varepsilon$. The Wigner function is defined in a quadratic manner, so it is insensitive to a constant phase shift.

The moments of the Wigner distribution function taken with respect to the momentum variable, provide the physical observables of the system. In particular, the position density and flux are given by
\begin{eqnarray}\label{eqn_Wigner-physObs}
 \rho^\varepsilon(\bd{x},t)&=|\psi^\varepsilon|^2=\int_{\R^d}F^\varepsilon \rd\bd{p},\quad J^\varepsilon(\bd{x},t)&=\varepsilon\text{Im}\left(\bar{\psi}^\varepsilon\cdot\nabla_{\bd{x}}\psi^\varepsilon\right)=\int_{\R^d}\bd{p}F^\varepsilon \rd\bd{p}.
\end{eqnarray}
The evolution of the Wigner function is governed by the Wigner equation
\begin{equation}\label{eqn_WignerEqn}
  \partial_t{F^\varepsilon}+\bd{p}\cdot\nabla_{\bd{x}}{F^\varepsilon}+\Xi[U\mbb{I}+V]F^\vep=0,
\end{equation}
where $\Xi[V]$ is defined as
\begin{align*}
 \begin{split}
  \Xi[V]F^\vep&=\f{1}{(2\pi)^d}\int_{\R^{2d}}\,\f{\ri}{\vep}\,\left[
  {V\left(\bd{x}-\f{\vep\bd{y}}{2}\right)F^\vep(\bd{x},\bd{p}')
     -F^\vep(\bd{x},\bd{p}')V\left(\bd{x}+\f{\vep\bd{y}}{2}\right)}\right]\e^{\ri(\bd{p}'-\bd{p})\cdot\bd{y}}\rd\bd{p}'\,\rd\bd{y}.
 \end{split}
\end{align*}
We note that $F$ and $V$ are matrices, and in general they do not commute.

A quantum mechanical operator can be univocally associated to a function $A(\bd{x},\bd{p})$ defined on the classical phase-space by the so called Weyl quantization \cite{weyl_quantenmechanik_1927,Morandi2010JPA}. The following map is used
\begin{equation}\label{def_WeylMap}
\mathcal{W}(A)[h](\bd{x})=\hat{A}[h](\bd{x})={\f{1}{(2\pi\vep)^d}}
\int\int\,A\left(\f{\bd{x}+\bd{y}}{2},\bd{p}\right)\,h(\bd{x},\bd{y})\,\e^{\f{\ri}{\vep}(\bd{x}-\bd{y})\cdot\bd{p}}\,\rd\bd{p}\,\rd\bd{y}.
\end{equation}
Here, $\hat{A}\equiv \mcl{W}(A)$ is the Weyl quantum mechanical operator defined on the space of the smooth functions $h(\bd{x},\bd{y})\in\mcl{S}(\R^d\times\R^d)$. The function $A(\bd{x},\bd{p})$ denotes the symbol of $\hat{A}$. It is easy to verify that the Weyl quantization map is the inverse of the Wigner transform (the Weyl quantization procedure applied to the Wigner function $F^\varepsilon$ provides the density operator). 

In particular, the Weyl quantization of the Hamiltonian is the Schr\"odinger operator.
Namely,
\begin{equation}\label{eqn_Hhat}
A(\bd{x},\bd{p}) = H(\bd{x},\bd{p})=\frac{\bd{p}^2}{2}+\tilde{V}(\bd{x})\quad\Rightarrow\quad
\hat{A} = \hat{H}=-\frac{\varepsilon^2}{2}\Delta_{\bd{x}}+\tilde{V}(\bd{x}).
\end{equation}
The use of the Wigner-Moyal formalism is eased by the definition of the Moyal product $\#$ as 
\begin{align}\label{def_MoyalProduct}
  A\#B:=&\frac{1}{(2\pi)^{2d}}\int{}A\left(\bd{x}-\frac{\varepsilon}{2}\boldsymbol\eta,\bd{p}+\frac{\varepsilon}{2}\boldsymbol\mu\right)B(\bd{x}',\bd{p}')e^{\ri(\bd{x}-\bd{x}')\cdot\boldsymbol\mu+\ri(\bd{p}-\bd{p}')\cdot\boldsymbol\eta}\rd\boldsymbol\mu\rd\bd{x}'\rd\boldsymbol\eta\rd\bd{p}' \nonumber \\
  =&A\e^{\f{\ri\vep}{2}\left(\overleftarrow\nabla_{\bd{x}}\cdot\overrightarrow\nabla_{\bd{p}}-\overleftarrow\nabla_{\bd{p}}\cdot\overrightarrow\nabla_{\bd{x}}\right)}B,
\end{align}
where the arrows indicate on which symbol the gradients act. An important property of the Moyal product is 
$\mathcal{W}(A\# B)=\mathcal{W}(A)\mathcal{W}(B)$. The \#-product admits an $\vep$-expansion. The $O(\vep)$ term is the classical Poisson bracket $\{A,B\}=\nabla_{\bd{p}}A\cdot\nabla_{\bd{x}}B-\nabla_{\bd{x}}A\cdot\nabla_{\bd{p}}B$, and
\begin{equation}\label{eqn_MoyalProduct_asym}
 A\#B=AB-\f{\ri\vep}{2}\{A,B\}+{O}(\vep^2).
\end{equation}

\subsection{The adiabatic case}\label{sec_Adiabatic}
The mathematical study of the semiclassical limit in the adaibatic case was carried out in \cite{lions_sur_1993,gerard_homogenization_1997}. 
According to Theorem 6.1 in \cite{gerard_homogenization_1997}, outside the
crossing set $S=\{\bd{x}:\,\ld^+(\bd{x})=\ld^-(\bd{x})\}$, the Wigner function can be obtained by the projection of the solution onto the eigenspaces of the Hamiltonian. 
Let $F^0(t,\bd{x},\bd{p})\doteq\lim_{\vep\rightarrow0}F(t,\bd{x},\bd{p})$, we have 
 \begin{subequations}
 \begin{align}
  F^0(t,\cdot)&=\Pi_+F^0(t,\cdot)\Pi_++\Pi_-F^0(t,\cdot)\Pi_-\label{eqn_decomposition}\\
  &=f^+(t,\cdot)\Pi_++f^-(t,\cdot)\Pi_-,\label{eqn_SimplifiedDecomposition}
 \end{align}
 \end{subequations}
where 
$$
\Pi_\pm(\bd{x})=\chi^{\pm}(\bd{x})\otimes\chi^{\pm}(\bd{x})
$$
and
$f^\pm$ are the particle densities related to energy levels $\ld^{\pm}(\bd{x},\bd{k})$
\begin{equation}
f^\pm=\mathrm{Tr}(\Pi_{\pm}F^0(t,\cdot)).
\end{equation}
The distributions $f^\pm$ satisfy  the classical Liouville equation:
\begin{subequations}\label{eqn_adiabatic}
 \begin{equation}
  \p_tf^\pm+\nabla_{\bd{p}}\ld^\pm\cdot\nabla_{\bd{x}}f^\pm-\nabla_{\bd{x}}\ld^\pm\cdot\nabla_{\bd{p}}f^\pm=0,\quad t>0,\;\bd{x}\in\R^d\setminus{S},\;\bd{p}\in\R^d,
 \end{equation}
 \begin{equation}
  f^\pm(t=0,\bd{x},\bd{p})=\mathrm{Tr}(\Pi_{\pm}F^0(t=0,\bd{x},\bd{p})).
 \end{equation}
\end{subequations}

\subsection{Quantum transition in the non-adiabatic case}\label{sec_Non-adiabatic}

In the proximity of the crossing points, the Born-Oppenheimer approximation is no longer valid. By using the Wigner formalism, we derive a semiclassical model that is able to treat the quantum mechanical band transitions in the case where the separation between the upper and lower energy levels scales as $\sqrt{\vep}$. In particular, the transitions between bands are captured by the off-diagonal terms of the Wigner matrix distribution. Our approach is alternative to the use of the  Landau-Zener formula for the evaluation of the transition probability in correspondence to an avoided crossing and overcomes some of the difficulties that affect the Landau-Zener approach as argued in~\cite{herman_generalization_1982}. We follow the derivation presented in~\cite{morandi_wigner_2011,kapral_mixed_1999}.

Starting from the von Neumann equation, we have:
\begin{align}\label{eqn_vonNeumannDiag}
  \ri\vep\pl{\hat{F}'}{t}&=\left[\hat{H}',\hat{F}'\right].
\end{align}
where we have defined
\begin{equation}
\hat{F}' = \hat{\Tht}\,\hat{\rho}\,\hat{\Tht}^\dagger,\quad \hat{H}'=\hat{\Tht}\,\hat{H}\,\hat{\Tht}^\dagger.
\end{equation}
Here, $\hat{\Tht}=\mcl{W}[\Tht]$ and $\hat{\Tht}^\dagger$ are the Weyl quantization of $\Tht$
and $\Tht^\dagger$ respectively. In particular, $\hat{\Tht}\hat{V}\hat{\Tht}^\dagger=\hat{\Lambda}_V$.

Equation~\eqref{eqn_vonNeumannDiag} is a diagonalized version of the von Neumann equation written on the Weyl operator formalism. To obtain an equivalent  dynamical system defined on the quantum phase-space, we usd the inverse Weyl mapping:
\begin{equation}\label{eqn_F'evolution}
 \ri\vep\pl{{F}'}{t}=\left[{H}',{F}'\right]_{\#},
\end{equation}
where $\left[A,B\right]_{\#}=A\#B-B\#A$ is a commutator of the Moyal
product~\eqref{def_MoyalProduct}, and $H'$ and $F'$ are symbols associated
with $\hat{H'}$ and $\hat{F'}$ respectively:
\begin{equation}\label{def_H'F'diagonal}
 H'(\bd{x},\bd{p})=\Tht(\bd{x})\#H(\bd{x},\bd{p})\#\Tht(\bd{x})^\dagger,\qquad F'(\bd{x},\bd{p})=\Tht(\bd{x})\#F(\bd{x},\bd{p})\#\Tht(\bd{x})^\dagger.
\end{equation}
By using~\eqref{eqn_MoyalProduct_asym}, we expand the
equation~\eqref{eqn_F'evolution}
\begin{align}\label{eqn_F'evolution_expand}
 \ri\vep\pl{{F}'}{t}
 =&[\Ld,F']-\f{\ri\vep}{2}\{\Ld,F\}+\f{\ri\vep}{2}\{F,\Ld\}
                      +\ri\vep[\bd{p}\cdot\nabla_{\bd{x}}\Tht\,\Tht^\dagger,F']\;+\;O(\vep^2)\nn\\
 \begin{split}
  =&[\Ld,F']-\f{\ri\vep}{2}[\nabla_{\bd{p}}\Ld,\nabla_{\bd{x}}F]_{+}
             +\f{\ri\vep}{2}[\nabla_{\bd{x}}\Ld,\nabla_{\bd{p}}F]_{+}\\
   &\quad\quad\quad        +\,\ri\vep\,[\bd{p}\cdot\nabla_{\bd{x}}\Tht\,\Tht^\dagger,F']\;+\;O(\vep^2),
 \end{split}
\end{align}
with $[A,B]_+=AB+BA$. Keeping up to the second order in $\vep$ (see Appendix for details):
\begin{equation}\label{eqn_H'diag}
 H'(\bd{x},\bd{p})=\Ld(\bd{x},\bd{p})+\ri\vep\bd{p}\cdot\nabla_{\bd{x}}\Tht(\bd{x})\Tht^\dagger(\bd{x})+\f{\vep^2}{2}\nabla_{\bd{x}}\Tht(\bd{x})\cdot\nabla_{\bd{x}}\Tht^\dagger(\bd{x}).
\end{equation}

By ignoring the $O(\varepsilon^2)$ terms, the evolution equations then become (for details of the asymptotic
derivation see Appendix)
\begin{subequations}\label{eqn_SemiclassicalLiouville_complex}
\begin{align}
 \pl{f^+}{t}&=
 -\bd{p}\cdot\nabla_{\bd{x}}f^++\nabla_{\bd{x}}\big(U+E\big)\cdot\nabla_{\bd{p}}f^+
 +\bar{b}^if^i+b^i\ov{f^i},\\
 \pl{f^-}{t}&=
 -\bd{p}\cdot\nabla_{\bd{x}}f^-+\nabla_{\bd{x}}\big(U-E\big)\cdot\nabla_{\bd{p}}f^-
 -\bar{b}^if^i-b^i\ov{f^i},\\
 \pl{f^i}{t}&=-\bd{p}\cdot\nabla_{\bd{x}}f^i\,+\nabla_{\bd{x}}U\cdot\nabla_{\bd{p}}f^i
 +b^i(f^--f^+)+(b^+-b^-)f^i+\f{2E}{\ri\vep}f^i,\label{eqn_SemiclassicalLiouville_complex_c}
\end{align}
\end{subequations}
where we have denoted
\begin{align}\label{def_F'components}
 F'&=
  \left(\begin{array}{cc}
         f^+ & f^i \\
         \bar{f}^i & f^-
        \end{array}\right),
 \quad\text{and}\quad
 \bd{p}\cdot\nabla_{\bd{x}}\Tht\,\Tht^\dagger=
  \left(\begin{array}{cc}
         b^+ & b^i \\
         -\bar{b}^i & b^-
        \end{array}\right).
\end{align}

In vector form \eqref{eqn_SemiclassicalLiouville_complex} becomes:
\begin{subequations}\label{eqn_SemiclassicalLiouville_complex_vector}
\begin{align}
 \pl{\bd{f}}{t}+\bd{p}\cdot\nabla_{\bd{x}}\bd{f}-\nabla_{\bd{x}}A\cdot\nabla_{\bd{p}}\bd{f}=C\bd{f}+\f{D}{\ri\vep}\bd{f}
\end{align}
\text{where:}
\begin{align}
 \bd{f}&=\left(f^+,\,f^-,\,f^i,\,\ov{f^i}\right)^T,\label{f_vect}\\
 A&=\text{diag}\left(U+E,\,U-E,\,U,\,U\right),\\
 D&=\text{diag}\left(0,\,0,\,2E,\,-2E\right),\\
 C&=\left(\begin{array}{cccc}
           0 & 0 & \ov{b^i} & b^i\\
           0 & 0 &-\ov{b^i} &-b^i\\
           -b^i & b^i & b^+-b^- & 0 \\
           -\ov{b^i} & \ov{b^i} & 0 & b^+-b^-
          \end{array}\right).\label{eqn_CouplingMatrix}
\end{align}
\end{subequations}
Here, $f^\pm$, both real, represent the projection coefficients onto the positive and
negative energy bands. The function $f^i$ describes the transition
between the two bands.

Specifically for the three examples
in~\eqref{eqn_u_v_delta_1d},~\eqref{eqn_u_v_delta_real} and~\eqref{eqn_u_v_delta}, we
have explicit formulae for $b^s$, $s\in\{\pm,i\}$:
\begin{description}
 \item[For~\eqref{eqn_u_v_delta_1d}:]
 $b^+=b^-\equiv0, \quad\text{and}\quad  b^i=-\frac{p\dt}{2(x^2+\delta^2)}$.
 \item[For~\eqref{eqn_u_v_delta_real}:] denote $\bd{p}=(p,q)$, then
 \begin{align}
  b^+=b^-\equiv0, \quad  b^i=\frac{1}{2(x^2+y^2+\delta^2)} \left( \frac{qxy}{\sqrt{y^2+\delta^2}} - p\sqrt{y^2+\delta^2} \right).
 \end{align}
 \item[For~\eqref{eqn_u_v_delta}:] one has for $\bd{p}=(p,q)$, $E=\sqrt{x^2+y^2+\dt^2}$ and
 \begin{align}
  \begin{split}
   b^+&=\frac{q \dt \left(E+x\right)}{2E^3}\ri\;,\quad
   b^-=\frac{q \dt}{2E(E+x)}\ri\;,\\
   b^i&=\frac{1}{2(x^2+y^2+\delta^2)} \left\{\left( \frac{qxy}{\sqrt{y^2+\delta^2}} - p\sqrt{y^2+\delta^2}\right) -\ri q\dt\right\}.
  \end{split}
 \end{align}
\end{description}

The system~\eqref{eqn_SemiclassicalLiouville_complex_vector} is hyperbolic, 
$\Theta$ is unitary, $b^\pm$ are purely imaginary, and the matrix $C$ is skew
Hermitian.

 \begin{remark}
 According to the adiabatic theory \cite{Niu2010Berry,morandi_wigner_2011}, when time $t$ is sufficiently small, the solution of the Schr\"odinger equation \eqref{eqn_Schrodinger} can be written as
  \begin{align}\label{eqn_state_in_one_band}
    &\psi^\vep(t)=\psi_+(t)+\psi_-(t),\nn\\
    \text{with }&\;\psi_\pm(t)=\e^{\ri\gamma_\pm(t)}\exp\left(\mp\f{\ri}{\varepsilon}\int_{0}^t\rd{t'}E(\bd{x}(t'))\right)\chi_\pm(\bd{x}(t)),
  \end{align}
  where $\bd{x}(t)$ is a semiclassical trajectory. For $t=0$, the initial state coincides with the eigenstate $\chi_s(\bd{x}(0))$ for $s\in\{+,-\}$. The second exponential in \eqref{eqn_state_in_one_band} is known as the dynamical phase factor, and $\gamma_\pm$ in the first exponential is the path integral of the Berry connection, i.e.
  \begin{align}
   \gamma_\pm(t)=\ri\int\,\dot{\bd{x}}(t)\cdot\big(\nabla_{\bd{x}}\chi_\pm(\bd{x}(t))\cdot\chi_\pm^\dagger(\bd{x}(t)\big)\,\rd{t},
  \end{align}
  which is called the Berry phase. This term cancels out in the diagonal term of the density function $\psi_+^\dagger\psi_+$ and $\psi_-^\dagger\psi_-$. However for the off diagonal term $\psi_{-}^\dagger{\psi}_+$, we have
 \begin{align}
  \begin{split}
   \psi_{-}^\dagger{\psi}_+(\bd{x}(t))=&\exp\left\{\ri\left(\gamma_+(t)-\gamma_-(t)-\f{2}{\varepsilon}\int_{0}^t\rd{t'}E(\bd{x}(t'))\right)\right\}\chi_-^\dagger(\bd{x}(t))\chi_+(\bd{x}(t))\,.
  \end{split}
  \end{align}
  By evaluating the derivative of the Berry phase we have
  \begin{align}
  \begin{split}
    &\ri\f{\rd}{\rd{t}}\left(\gamma_+(t)-\gamma_-(t)-\f{2}{\varepsilon}\int_{0}^t\rd{t'}E(\bd{x}(t'))\right)\\*
   =&-\dot{\bd{x}}(t)\cdot\big(\nabla_{\bd{x}}\chi_+\cdot\chi_+^\dagger-\nabla_{\bd{x}}\chi_-\cdot\chi_-^\dagger\big)-\f{2\ri}{\varepsilon}E(\bd{x}(t)).
  \end{split}
  \end{align}
  If we apply $\dot{\bd{x}}=\bd{p}$ where $\bd{p}$ is the momentum, then we get
  \begin{align}
  \begin{split}
    &\ri\f{\rd}{\rd{t}}\left(\gamma_+(t)-\gamma_-(t)-\f{2}{\varepsilon}\int_{0}^t\rd{t'}E(\bd{x}(t'))\right)
   =\big(b^+(\bd{x}(t))-b^-(\bd{x}(t))\big)-\f{2\ri}{\varepsilon}E(\bd{x}(t))\,.
  \end{split}
  \end{align}
  Comparing with \eqref{eqn_SemiclassicalLiouville_complex_c}, one can see that these are {\it exactly} the coefficients of the $f^i$ terms in \eqref{eqn_SemiclassicalLiouville_complex_c}. This shows our model indeed captures the Berry phase in the inter-band transition processes.
 
 \end{remark}

\section{A hybrid model by domain decomposition}\label{sec_NumDDM}
The equation~\eqref{eqn_SemiclassicalLiouville_complex} is a hyperbolic system, with a transport part and a source term. Concerning the numerical treatment of Eq. \eqref{eqn_SemiclassicalLiouville_complex} the major difficulties arise form the term  $\frac{2E}{\ri\varepsilon}$ in the equation for $f^i$. It introduces rapid oscillations in both space and time that demand high computational cost. In order to reduce the numerical complexity, we solve the semiclassical model \eqref{eqn_SemiclassicalLiouville_complex} only in the proximity of the crossing zone. Away form the crossing points we neglect the band transitions and solve the adiabatic model
\begin{subequations}\label{eqn_SemiclassicalLiouville_adiabatic_xy}
 \begin{eqnarray}
  \pl{f^+}{t}&=&-\bd{p}\cdot\nabla_{\bd{x}}f^++\nabla_{\bd{x}}\big(U+E\big)\cdot\nabla_{\bd{p}}f^+,\\
  \pl{f^-}{t}&=&-\bd{p}\cdot\nabla_{\bd{x}}f^-+\nabla_{\bd{x}}\big(U-E\big)\cdot\nabla_{\bd{p}}f^-.
 \end{eqnarray}
\end{subequations}
A similar hybrid model was used in~\cite{chai_semi-classical_2013} for the Schr\"odinger equation with a periodic lattice potential.

As an example, consider the one-dimensional case with $p>0$ (so that in the $x$ space the wave packet moves from the left to the right). The other cases are treated similarly. We decompose the domain into the following two regions: 

 \begin{description}
  \item[The adiabatic region:] $x<-C_0\sqrt{\varepsilon}$ and $x>C_0\sqrt{\varepsilon}$:
  \newline
          In this region, we use $o(1)$ coarse mesh, independent of
          $\varepsilon$ for the adiabatic Liouville system
          \eqref{eqn_SemiclassicalLiouville_adiabatic_xy}. $f^i$ is set to be
          zero. At $x=-C_0\sqrt{\varepsilon}$, no boundary condition is required, while at $x=C_0\sqrt{\varepsilon}$,  we impose the inflow boundary condition that
          $f^+$ and $f^-$ are given by the solution inside the non-adiabatic region discussed below.
  \item[The non-adiabatic region:]
  $\left[-C_0\sqrt{\varepsilon},\quad C_0\sqrt{\varepsilon}\right]$:
  \newline
          In this region we use $o(\sqrt{\varepsilon})$ mesh and compute the
          full system~\eqref{eqn_SemiclassicalLiouville_complex}.  The system is hyperbolic, so
          the boundary condition only needs to be specified in the incoming direction. The incoming
          boundary data for $f^i$ is set to be zero.
          Since the region size is of
          $\mathcal{O}(\sqrt{\varepsilon})$, the total number of grid points along $x$-direction
          remain independent from $\varepsilon$.
\end{description}
In our simulation we choose $C_0 = 3$.

\section{Numerical examples}
In our numerical simulation, we use the hybrid model proposed in previous section. For the transport operator of the Liouville systems (both adiabatic and non-adiabatic cases), we use the standard second order upwind total-variation-diminishing (TVD) scheme with van Leer slope limiter \cite{LeV}. The reference solutions are obtained by the direct computation of the Schr\"odinger equation~\eqref{eqn_Schrodinger} with the time-splitting spectral method described in~\cite{bao_time-splitting_2002}. In the examples presented in the next sections, we use the following initial data:
\begin{equation}\label{eqn_initial_Sch_bands}
 \psi^\vep(t=0,\bd{x})=\psi^\vep_0(\bd{x})=g^\vep_0(\bd{x})\big(a^+\chi^+(\bd{x})+a^-\chi^-(\bd{x})\big),
\end{equation}
where $g^\vep_0$ is the $\vep$-scaled Gaussian packet:
\begin{equation}\label{eqn_initial_Sch_packet}
 g^\vep_0(\bd{x})=\left(\f{A}{\pi}\right)^{d/4}\exp\left\{-\f{A}{2}|\bd{x}-\bd{x_{0}}|^2+\f{\ri}{\vep}\,\bd{p_{0}}\cdot(\bd{x}-\bd{x_{0}})\right\}.
\end{equation}
Here, $a^\pm$ are constants and $\chi^\pm$ the eigenvectors of the operator $\hat{V}$ (see Eq. \eqref{eqn_potential}). By using the definition of the Wigner transform~\eqref{def_WignerTransform} and \eqref{def_H'F'diagonal}, we obtain the initial condition for $F'$. In the regime  $\varepsilon\ll{1}$ we obtain
\begin{subequations}
\begin{align}\label{eqn_initial_fs}
 f^+(t=0,\bd{x},\bd{p})&={(a^+)^2}\left(\f{A}{\pi^2\vep}\right)^{d/2}\exp\left\{-A|\bd{x}-\bd{x_{0}}|^2-\f{1}{\vep}|\bd{p}-\bd{p_{0}}|^2\right\},\\
 f^-(t=0,\bd{x},\bd{p})&={(a^-)^2}\left(\f{A}{\pi^2\vep}\right)^{d/2}\exp\left\{-A|\bd{x}-\bd{x_{0}}|^2-\f{1}{\vep}|\bd{p}-\bd{p_{0}}|^2\right\},\\
 f^i(t=0,\bd{x},\bd{p})&=0.
\end{align}
\end{subequations}
where, according to \eqref{f_vect} we have expressed the initial data in terms of the components of the vector $\bd{f}$. In particular, we note that in the limit $\varepsilon \to 0$,  $f^+$ and $f^-$ become the classical Dirac measure $\delta(\bd{p}-\bd{p}_0)$.

In our numerical experiments, the relevant physical observables are the particle density in the lower ($-$) and upper ($+$) bands. In order to compare the solution of our  new model with the original Schr\"odinger equation, it is convenient to consider the expression of the particle density  in the two formulations
\begin{align}\label{def:rho}
\begin{cases}
 \rho^\pm_{schr}(t,\bd{x})&=|\Pi_\pm\psi^\vep(t,\bd{x})|^2, \quad\text{and}\quad P^\pm_{schr}=\int_{\R^d}\rho^\pm_{schr}(t,\bd{x})\,\rd\bd{x},\\
 \rho^\pm_{liou}(t,\bd{x})&=\int_{\R^d}f^\pm(t,\bd{x},\bd{k})\rd\bd{k}, \quad\text{and}\quad P^\pm_{liou}=\int_{\R^d}\rho^\pm_{liou}(t,\bd{x})\,\rd\bd{x}.
 \end{cases}
\end{align}
The total density is given by ~\cite{gosse_multiphase_2004,jin_computational_2008}:
\begin{equation}\label{def_Cumulative}
 M_{schr} = \int_{\Omega_x} \big(\rho^+_{schr}(\bd{y})+\rho^-_{schr}(\bd{y})\big)\rd\bd{y},\quad
 M_{liou} = \int_{\Omega_x}\big(\rho^+_{liou}(\bd{y})+\rho^-_{liou}(\bd{y})\big)\rd\bd{y},
\end{equation}
where $\Omega_{\bd{x}}=\left\{\bd{y}=(y_1,\cdots,y_d)\in\R^d\:\;:\;y_i\leq x_i,\;i=1,\cdots,d\right\}$
for $\bd{x}=(x_1,\cdots,x_d)\in\R^d$, and $\Omega\subset\R^d$ is the computational domain, $|\Omega|$ is the
measure of $\Omega$. In order to estimate the accuracy of our method, we define the following $L^1-$norm of the error
\begin{equation}\label{def_CumulativeErr}
 \text{Err}^\varepsilon=\f{1}{|\Omega|}\int_{\Omega} \left|M_{schr} - M_{liou}\right|\rd\bd{x},
\end{equation}

\subsection{One dimensional system}
In our first example, we consider a one dimensional system. We choose $u$ and $v$ as
\begin{equation}\label{eqn_example_uv}
 u(x)=x,\quad\text{and}\quad v(x)\equiv\delta=\f{\sqrt\vep}{4}.
\end{equation}
The minimum of the energy gap is $2\delta$ ($x=0$). The initial data for the Schr\"odinger equation are given in
\eqref{eqn_initial_Sch_bands}-\eqref{eqn_initial_Sch_packet} with $a^+=1$,
$a^-=0$, $x_0=0.3125$, and $p_0=-1$ (pure state initial condition). For the Schr\"odinger equation we use a uniform grid for the space and the time variables with, respectively, $\Delta x=\vep/32$ and $\Delta t=\vep/32$. For the semiclassical Liouville system~\eqref{eqn_SemiclassicalLiouville_complex}, the phase-space $(x,p)$ domain is discretized with a uniform mesh with $\Delta x = \Delta p = 2^{-9}$ in the adiabatic region  ($x< -3\sqrt\vep \cup x>  3\sqrt\vep$), and $\Delta x = \Delta p = 2^{-11}$ in the non-adiabatic region. The time step is chosen as $\Delta t=2^{-14}$.

\begin{figure}[htp]
  \begin{center}
   \includegraphics[scale = 0.8]{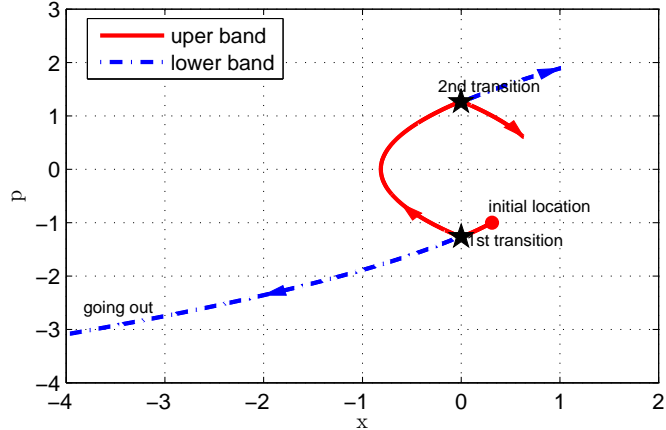}
  \end{center}
  \caption{The trajectories of wave packets: initially, the wave packet centers around the bullet point and starts moving towards the crossing point at $x=0$, marked as a star. It splits into two parts there: one of them, denoted by the dash-point line, keeps moving towards the left, while the other wave packet, the one that jumps up to the higher energy band, bounces back and hits the origin at $x=0$. Over there, it goes through the ``second transition'' and splits up into trajectories.}
  \label{fig:Trajectory}
 \end{figure}

\begin{figure}[htp]
  \begin{center}
   \includegraphics[width=4.5in]{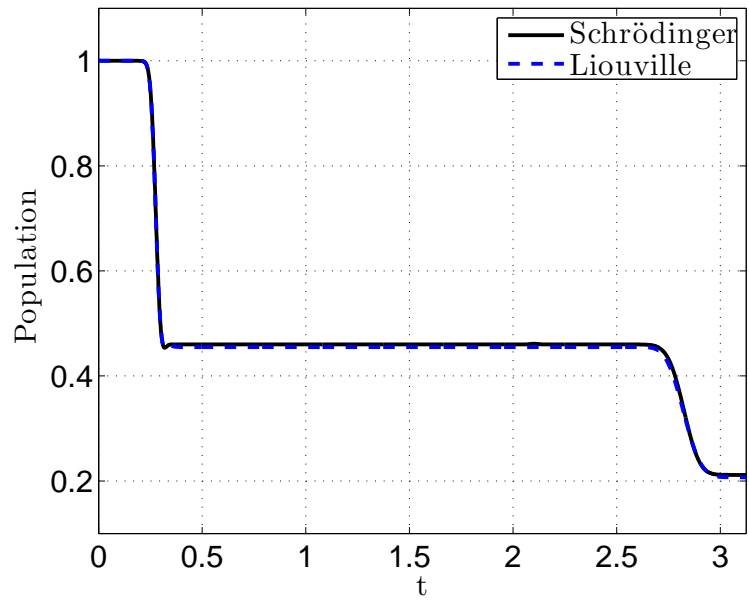}
  \end{center}
  \caption{Time evolution of $P^+$~\eqref{def:rho}.}
  \label{fig:D1Pure_longtime}
 \end{figure}

The evolution of the particle wave packets can be easily understood. The trajectories of the wave packets center-of-mass are illustrated in Figure~\ref{fig:Trajectory}. The Gaussian profile has a negative mean  velocity and passes through the crossing point $x=0$ at around $t = 0.25$. The Gaussian wave function splits into two parts. Around one half of the particles enter into the lower energy band, and the others stay in the upper band. The packet on the lower energy level is accelerated and leaves the simulation domain. The particles on the higher energy band are accelerated on the opposite direction, the momentum decreases and the wave packet is reflected around $x=-1$. At $t = 2.75$, the wave passes through the crossing point for the second time, and undergoes another hopping process. In Figure~\ref{fig:D1Pure_longtime} we show the evolution of $P^+$ with respect to time. As the wave packet passes through the crossing point twice, the mass gets transferred to another energy band twice, generating two 
jumps in $P^+$. The numerical results given by the semiclassical model shows good agreement with that of the Schr\"odingier equation.

\begin{figure}[htp]
  \begin{center}
   \subfloat[$\rho^+_{schr/liou}$]
   {\includegraphics[width=5in]{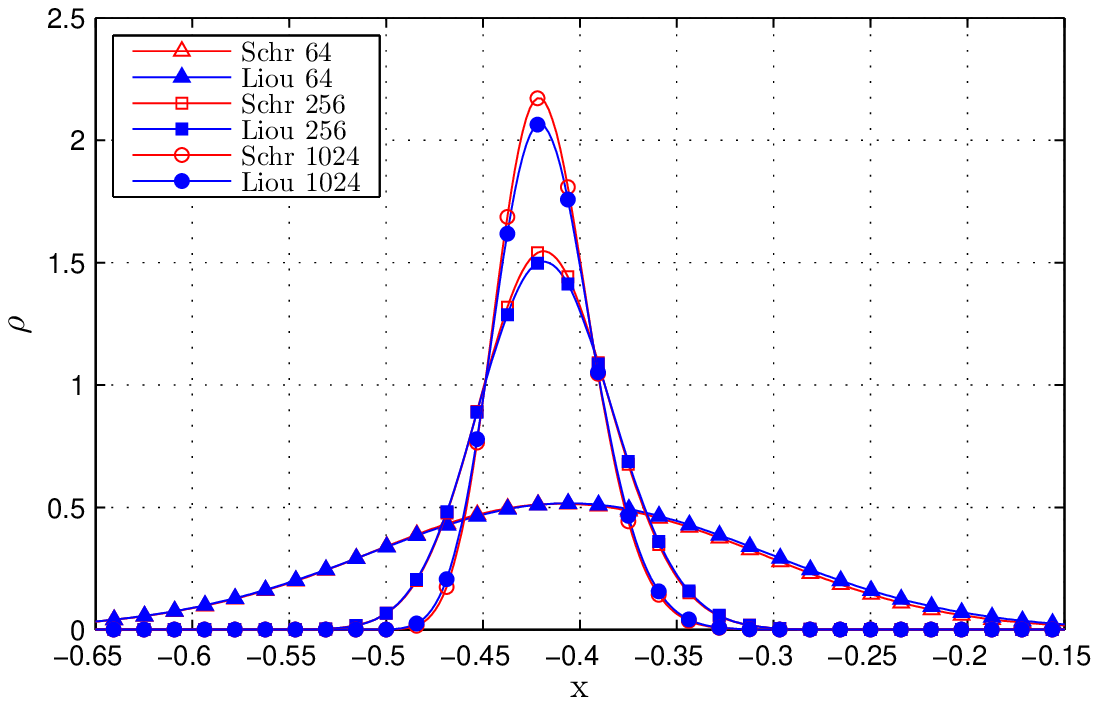}}\\
   \subfloat[$\rho^-_{schr/liou}$]
   {\includegraphics[width=5in]{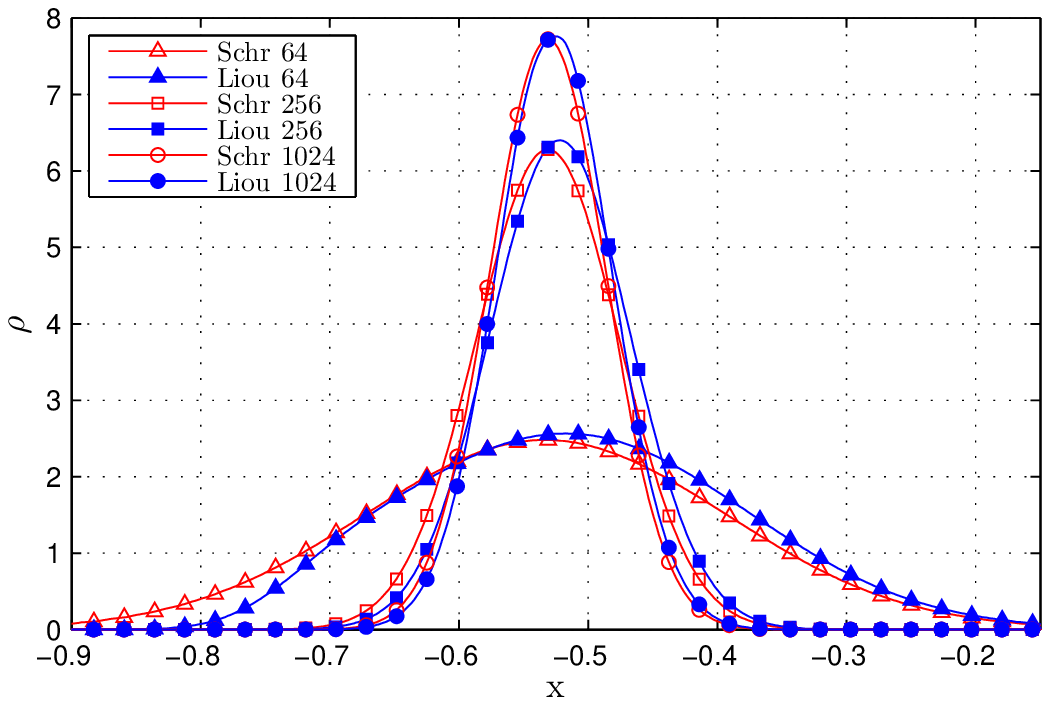}}
  \end{center}
  \caption{One dimensional system simulation: the density functions $\rho^{\pm}_{schr/liou}$ for different $\varepsilon$ at time $t=0.75$, $\delta=\sqrt\vep/4$. The legend ``Schr $n$'' (or ``Liou $n$'') represents the solution of the Schr\"odinger equation (or the hybrid model) with $\vep=1/n$. Here, $a^+=1$, $a^-=0$, $x_0=0.5$, and $p_0=-1$.}
  \label{fig:D1Pure_density}
 \end{figure}

\begin{figure}[htp]
  \begin{center}
    {\includegraphics[width=3.5in]{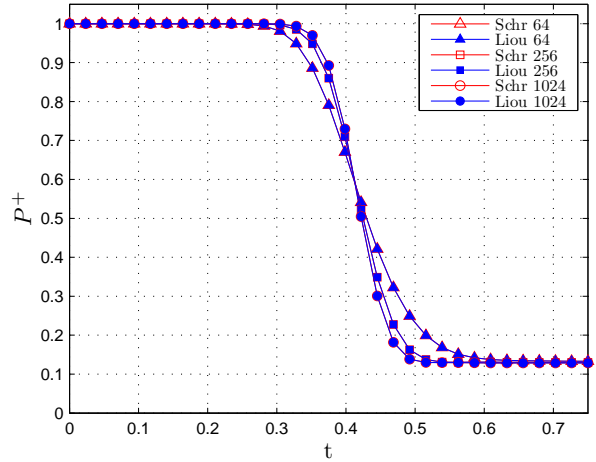}}
  \end{center}
  \caption{One dimensional system simulation: time evolution of the population on the upper band
  $P^+_{schr/liou}$ with $\delta=\sqrt\vep/4$. The legend ``Schr $n$'' (or ``Liou $n$'') represents the solution of the Sch\"odinger equation (or the hybrid model) with $\vep=1/n$. Here, $a^+=1$, $a^-=0$, $x_0=0.5$, and $p_0=-1$.}
  \label{fig:D1Pure_mass}
 \end{figure}

\begin{figure}[htp]
  \begin{center}
   \includegraphics[width=3.5in]{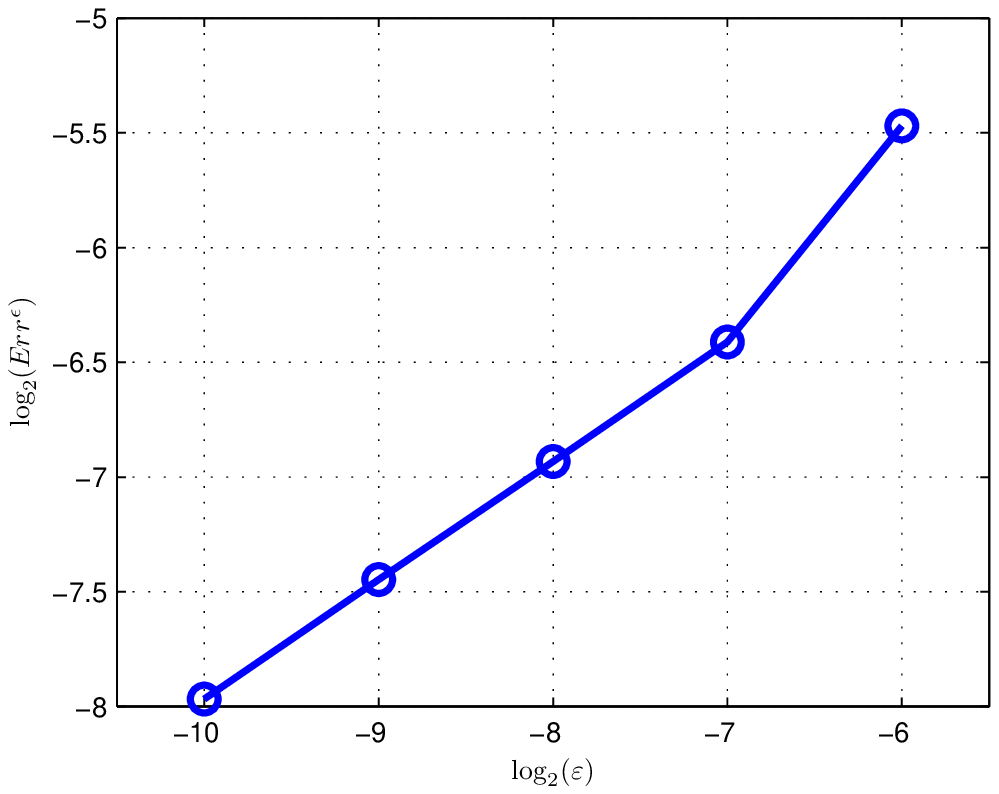}
  \end{center}
  \caption{One dimensional system simulation: $\text{Err}^\varepsilon$~\eqref{def_CumulativeErr} decreases with a rate of $O(\sqrt\varepsilon)$. Here, $a^+=1$, $a^-=0$, $x_0=0.5$, $p_0=-1$, and $t=0.75$.}
  \label{fig:D1Pure_err}
 \end{figure}

We compare the results of our model with that of the Schr\"odinger. In Figure~\ref{fig:D1Pure_density} we compare the results $\rho^\pm$ given by the two systems at time $t=0.75$ ($a^+=1$, $a^-=0$, $x_0=0.5$, and $p_0=-1$). In Figure~\ref{fig:D1Pure_mass} we check the evolution of the population on the first band $P^+$ along the time. Figure~\ref{fig:D1Pure_err} shows that the hybrid model error \eqref{def_CumulativeErr} decreases at the rate of $\mathcal{O}(\sqrt\varepsilon)$ for $\delta=\mathcal{O}(\sqrt{\varepsilon})$. The simulations show a good agreement between the two system.

We consider now a different initial condition. The initial wave packet for the Schr\"odinger equation is now given by~\eqref{eqn_initial_Sch_bands}-\eqref{eqn_initial_Sch_packet} with $a^+=a^-=1/\sqrt{2}$, $x_0=0.5$, and $p_0=-1$. This initial datum corresponds to a linear superposition of two Gaussian packets that belong to the upper and lower bands respectively.

\begin{figure}[htp]
  \begin{center}
   \subfloat[$\rho^+_{schr/liou}$]
   {\includegraphics[width=5in]{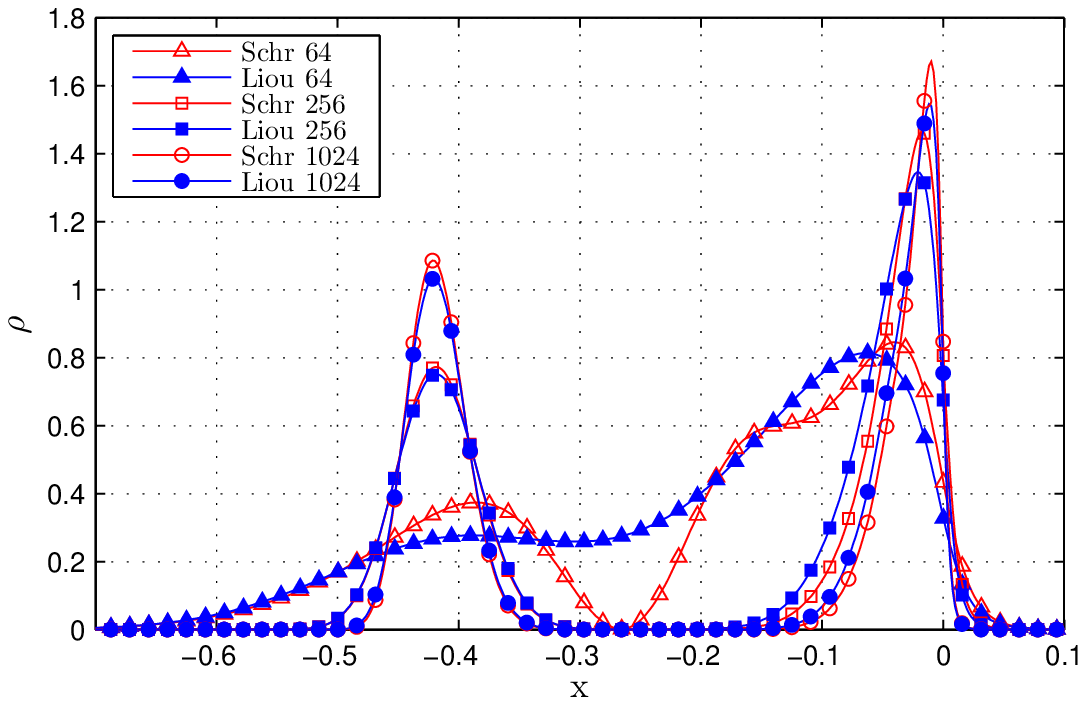}}\\
   \subfloat[$\rho^-_{schr/liou}$]
   {\includegraphics[width=5in]{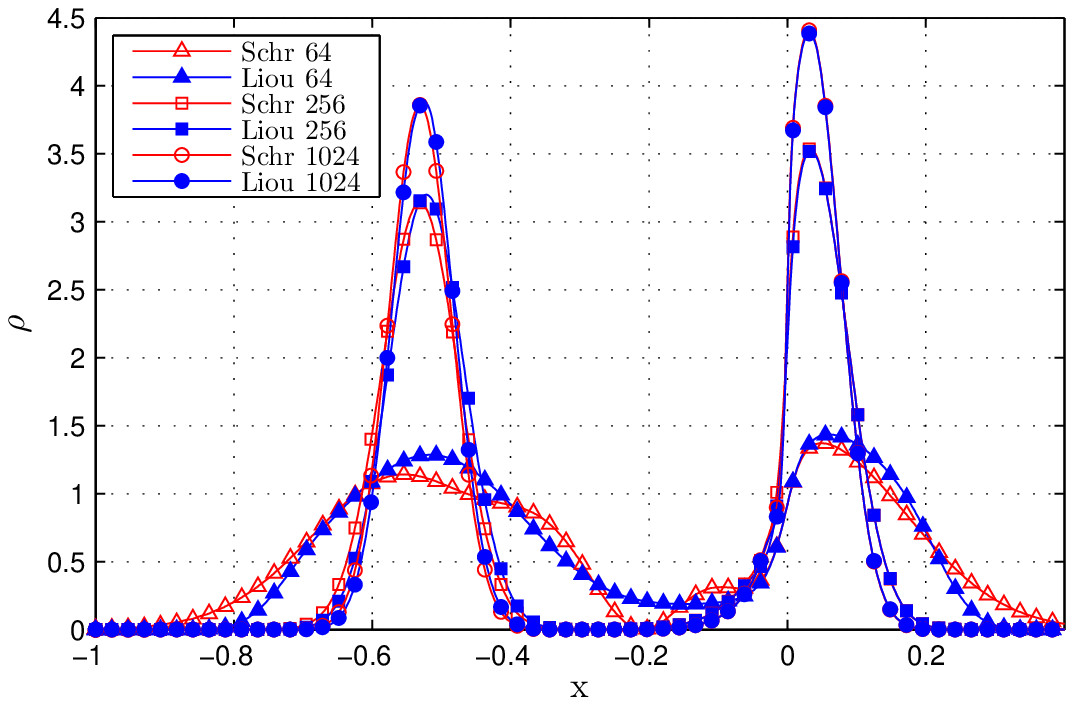}}
  \end{center}
  \caption{One dimensional system simulation: the density functions $\rho^{\pm}_{schr/liou}$ for different $\varepsilon$ at time $t=0.75$, $\delta=\sqrt\vep/4$. The legend ``Schr $n$'' (or ``Liou $n$'') represents the solution of the Sch\"odinger equation (or the hybrid model) with $\vep=1/n$. Here, $a^+=a^-=1/\sqrt{2}$, $x_0=0.5$, and $p_0=-1$.}
  \label{fig:D1Mix_density}
 \end{figure}

\begin{figure}[htp]
  \begin{center}
    \subfloat[$\varepsilon=2^{-6}$]
    {\includegraphics[width=2.5in]{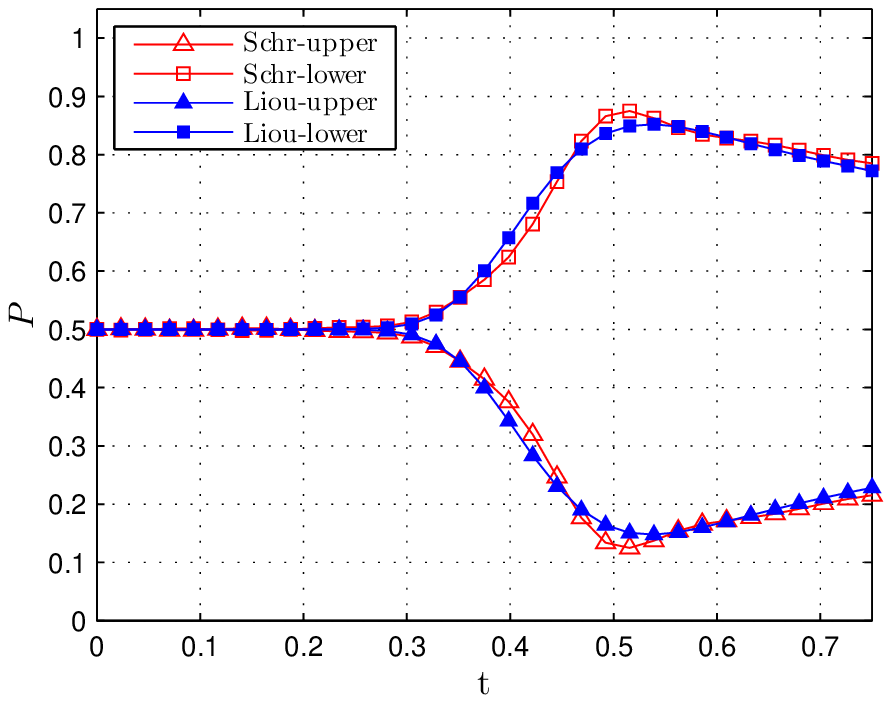}}
    \subfloat[$\varepsilon=2^{-7}$]
    {\includegraphics[width=2.5in]{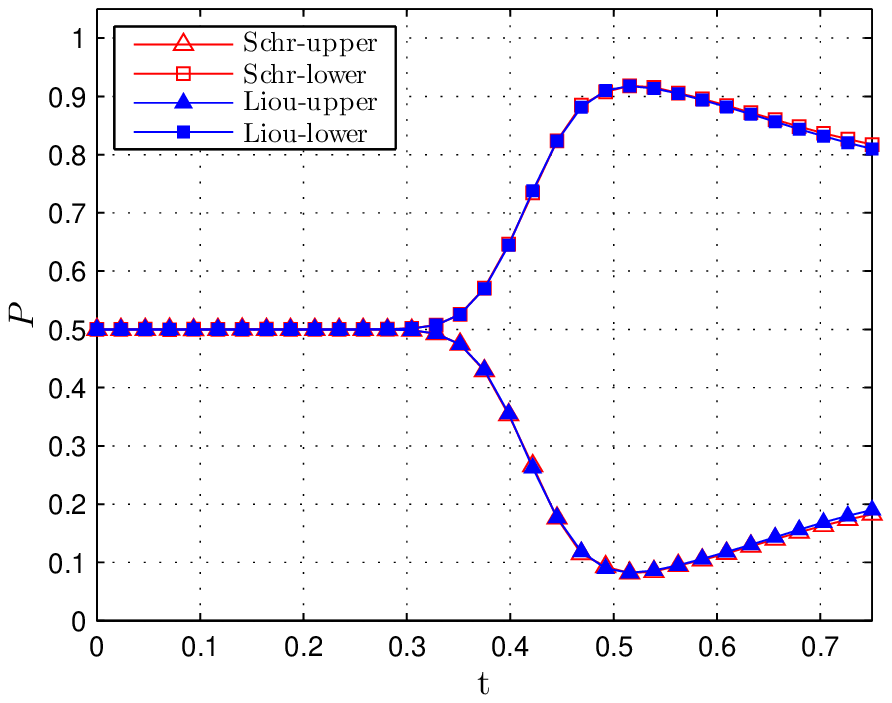}}\\
    \subfloat[$\varepsilon=2^{-8}$]
    {\includegraphics[width=2.5in]{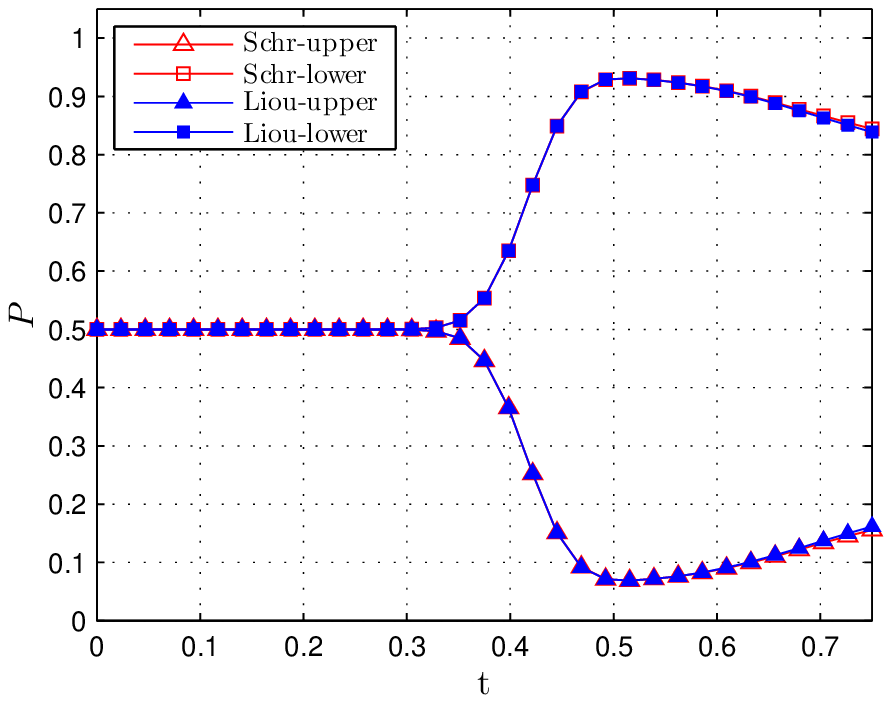}}
    \subfloat[$\varepsilon=2^{-9}$]
    {\includegraphics[width=2.5in]{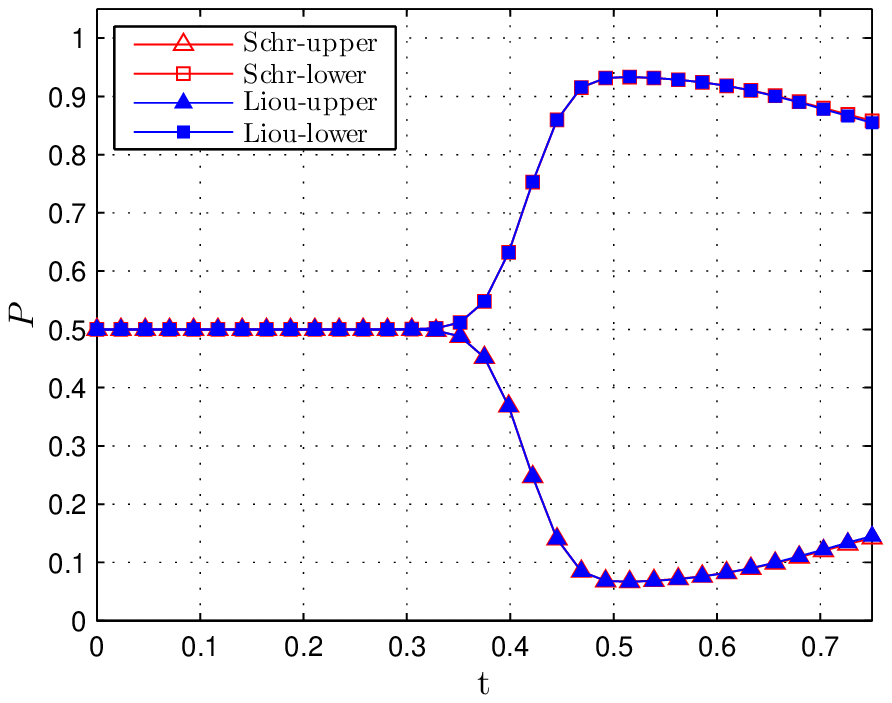}}
  \end{center}
  \caption{One dimensional system simulation: time evolution of the population on the upper and lower bands
  $P^\pm_{schr/liou}$~\eqref{def:rho}. $\delta=\sqrt\vep/4$. The legend
  ``Schr-upper'' ( or ``Schr-lower'') represents the population on upper (or lower) band given by the Schr\"odinger equation,
  ``Liou-upper'' ( or ``Liou-lower'') represents the population on upper (or lower) band given by the hybrid model. Here, $a^+=a^-=1/\sqrt{2}$, $x_0=0.5$, and $p_0=-1$.}
  \label{fig:D1Mix_mass}
 \end{figure}

\begin{figure}[htp]
  \begin{center}
   \includegraphics[width=3.5in]{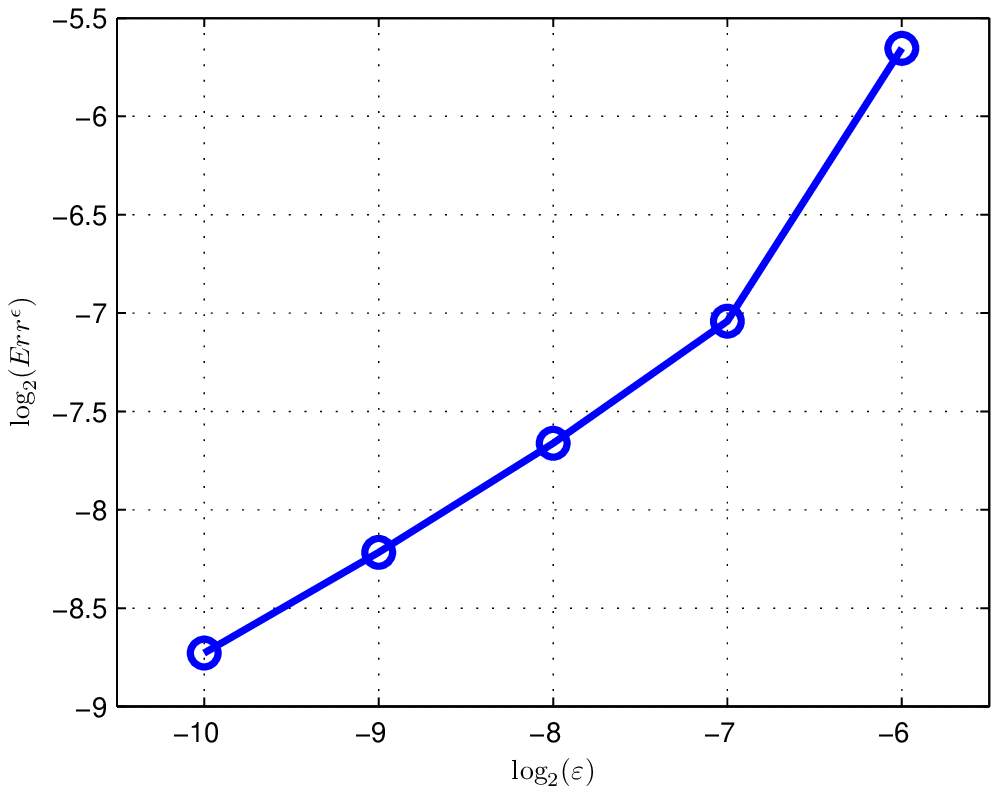}
  \end{center}
  \caption{One dimensional system simulation: $\text{Err}^\varepsilon$~\eqref{def_CumulativeErr} as a function of $\varepsilon$
  at $t=0.75$. Here, $a^+=a^-=1/\sqrt{2}$, $x_0=0.5$, and $p_0=-1$.}
  \label{fig:D1Mix_err}
 \end{figure}

In Figure~\ref{fig:D1Mix_density} we compare the numerical results of $\rho^\pm$ to the Schr\"odinger equation and those to~\eqref{eqn_SemiclassicalLiouville_complex} for $t=0.75$. In Figure~\ref{fig:D1Mix_mass}, we show the evolution of the populations on the upper and lower bands $P^\pm$ with respect to time. We see that in the case of large $\varepsilon$ the semiclassical solution is not completely satisfactory. However it is able to capture the main structure of the quantum interference between the upper and the lower band waves. For small $\vep$, the two wave packets are well-separated, the solution of the hybrid model is in good agreement with the Schr\"odinger solution. Figure~\ref{fig:D1Mix_err} shows that the cumulative error decreases at the rate of $\mathcal{O}(\sqrt{\varepsilon})$.

\subsection{Two-dimensional  system}
\begin{figure}[htp]
  \begin{center}
   \subfloat[$t=0.0781$]
    {\includegraphics[width=2.5in]{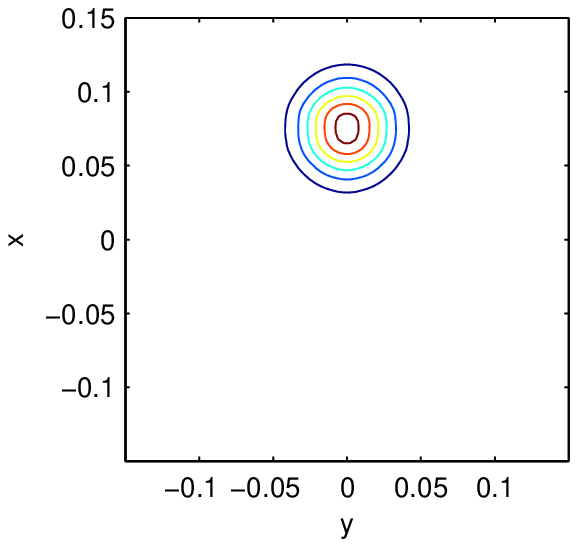}
     \includegraphics[width=2.5in]{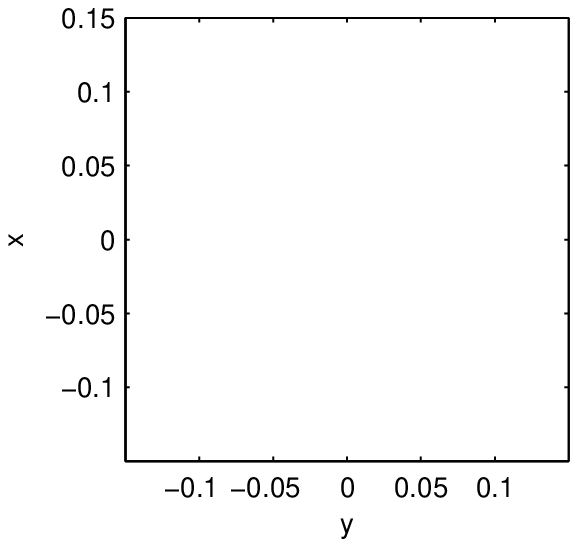}}\\
   \subfloat[$t=0.1562$]
    {\includegraphics[width=2.5in]{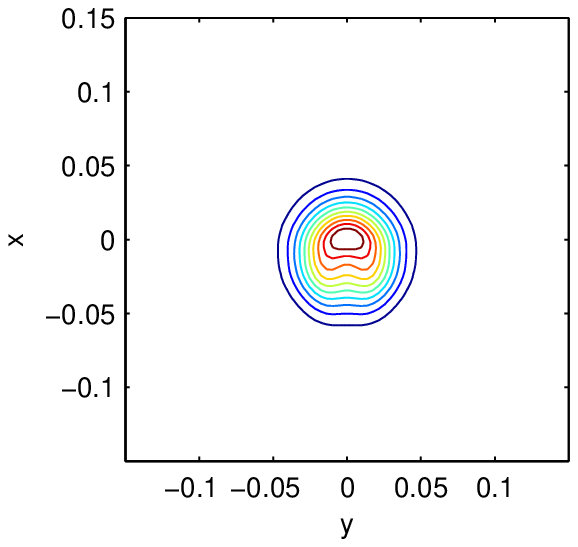}
     \includegraphics[width=2.5in]{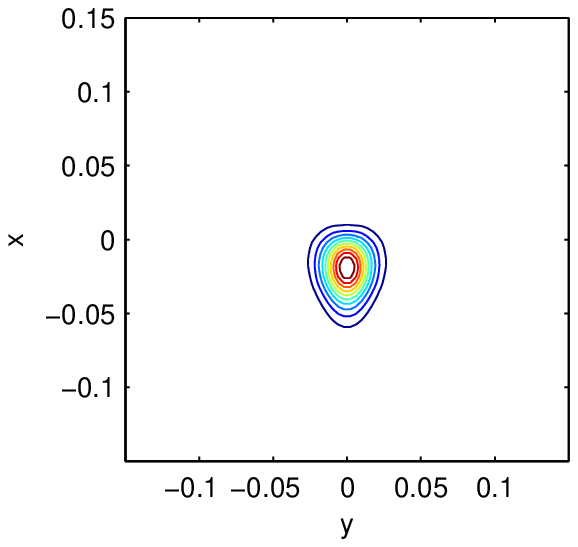}}\\
   \subfloat[$t=0.2344$]
    {\includegraphics[width=2.5in]{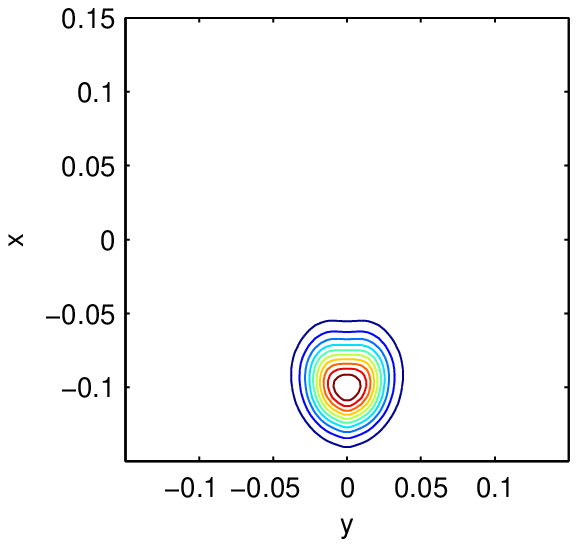}
     \includegraphics[width=2.5in]{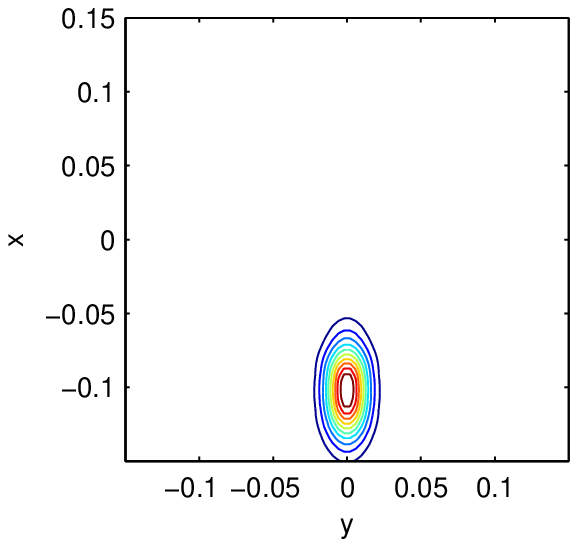}}
  \end{center}
  \caption{Two dimensional system simulation: the time evolution of density contour computed by the hybrid model. The left/right column are for $\rho^+_{liou}$/
  $\rho^-_{liou}$, the density on the upper/lower band~\eqref{def:rho}.
  $\delta=\sqrt\vep/2$ and $\vep=2^{-10}$. One can see in (b) around time $t=0.1562$, the wave packet hits the crossing point and a portion of the mass jumps to the upper band.}\label{fig:D2Pure_density}
 \end{figure}

\begin{figure}[htp]
  \begin{center}
    \subfloat[$\varepsilon=2^{-6}$]
    {\includegraphics[width=2.5in]{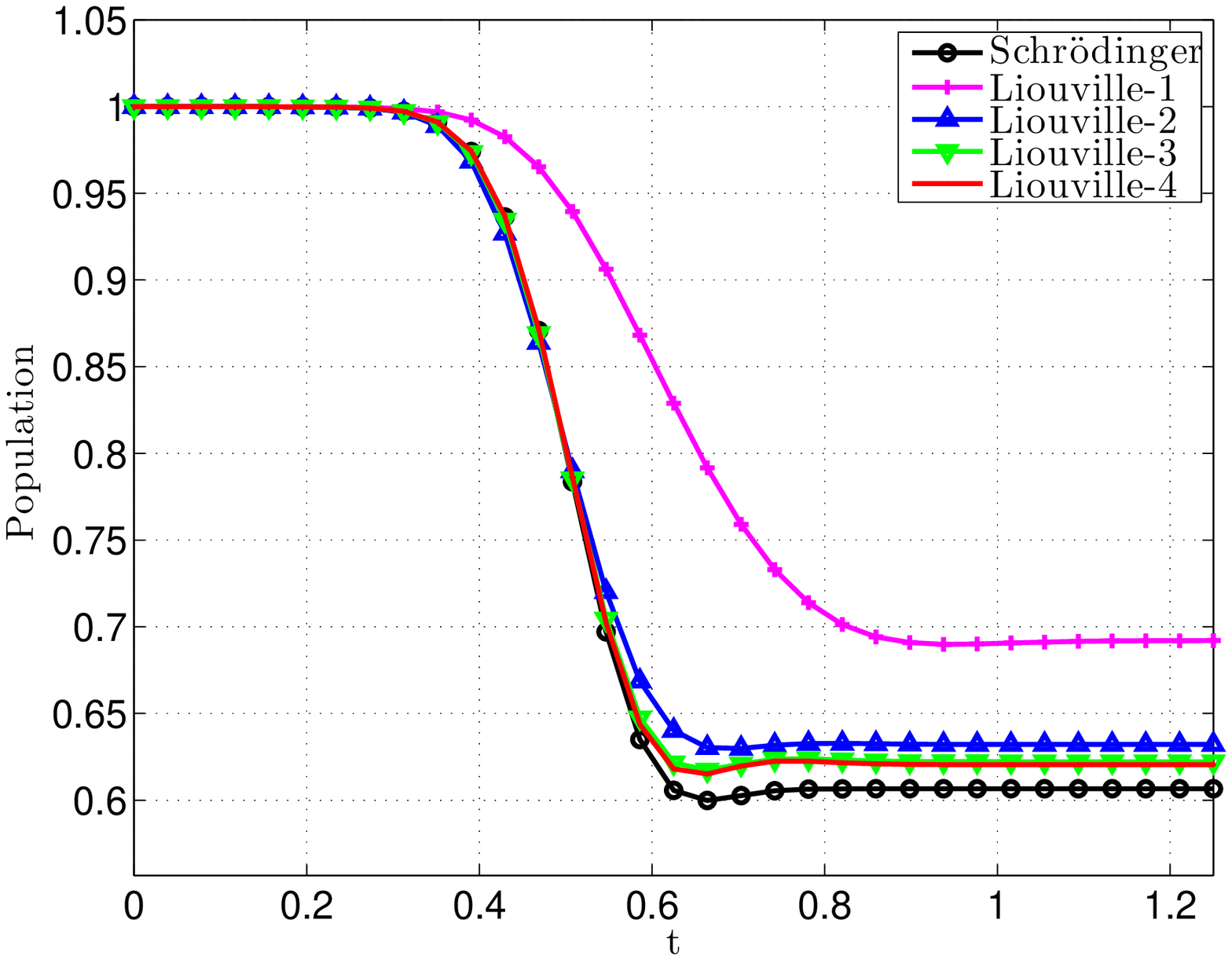}}
    \subfloat[$\varepsilon=2^{-7}$]
    {\includegraphics[width=2.5in]{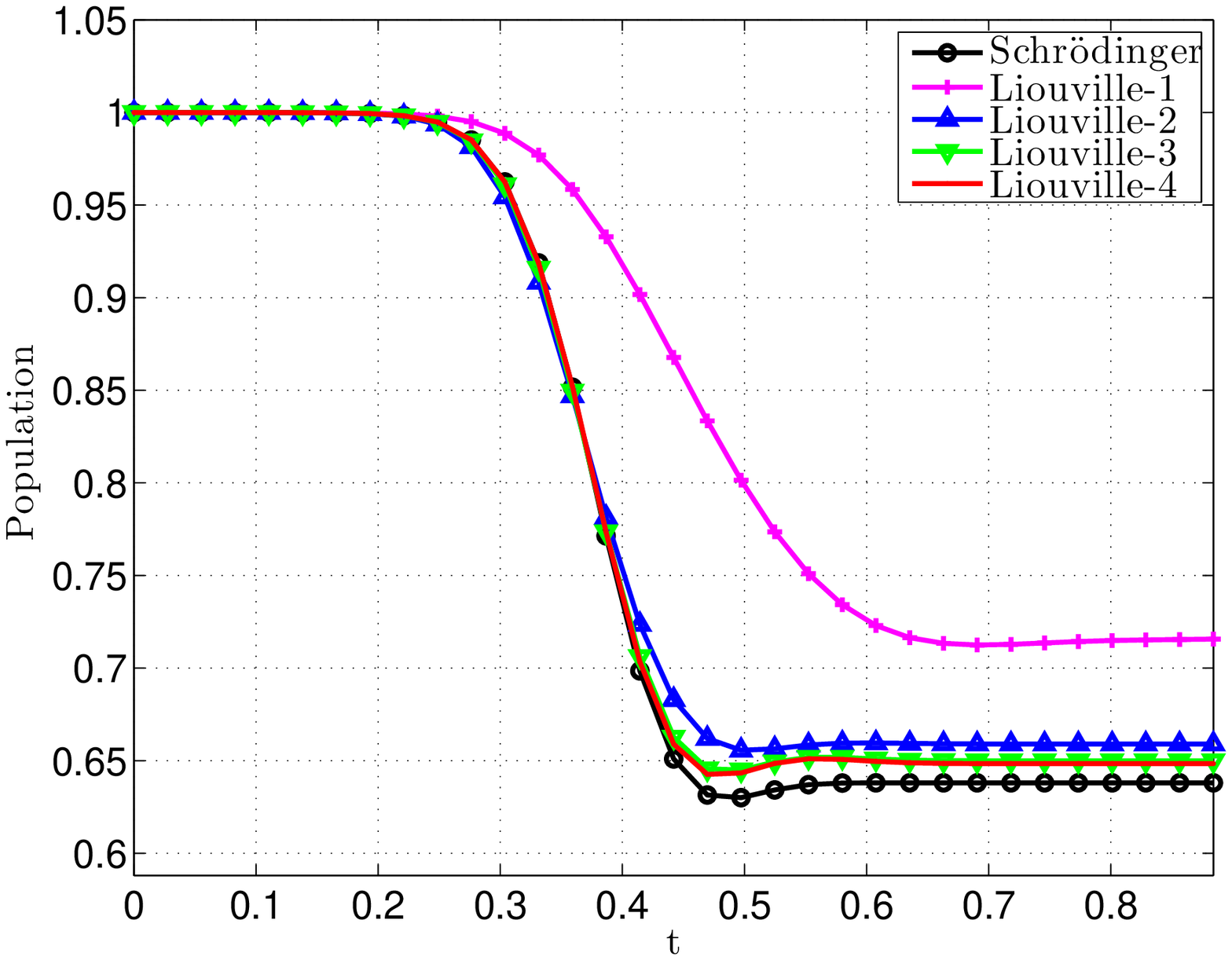}}\\
    \subfloat[$\varepsilon=2^{-8}$]
    {\includegraphics[width=2.5in]{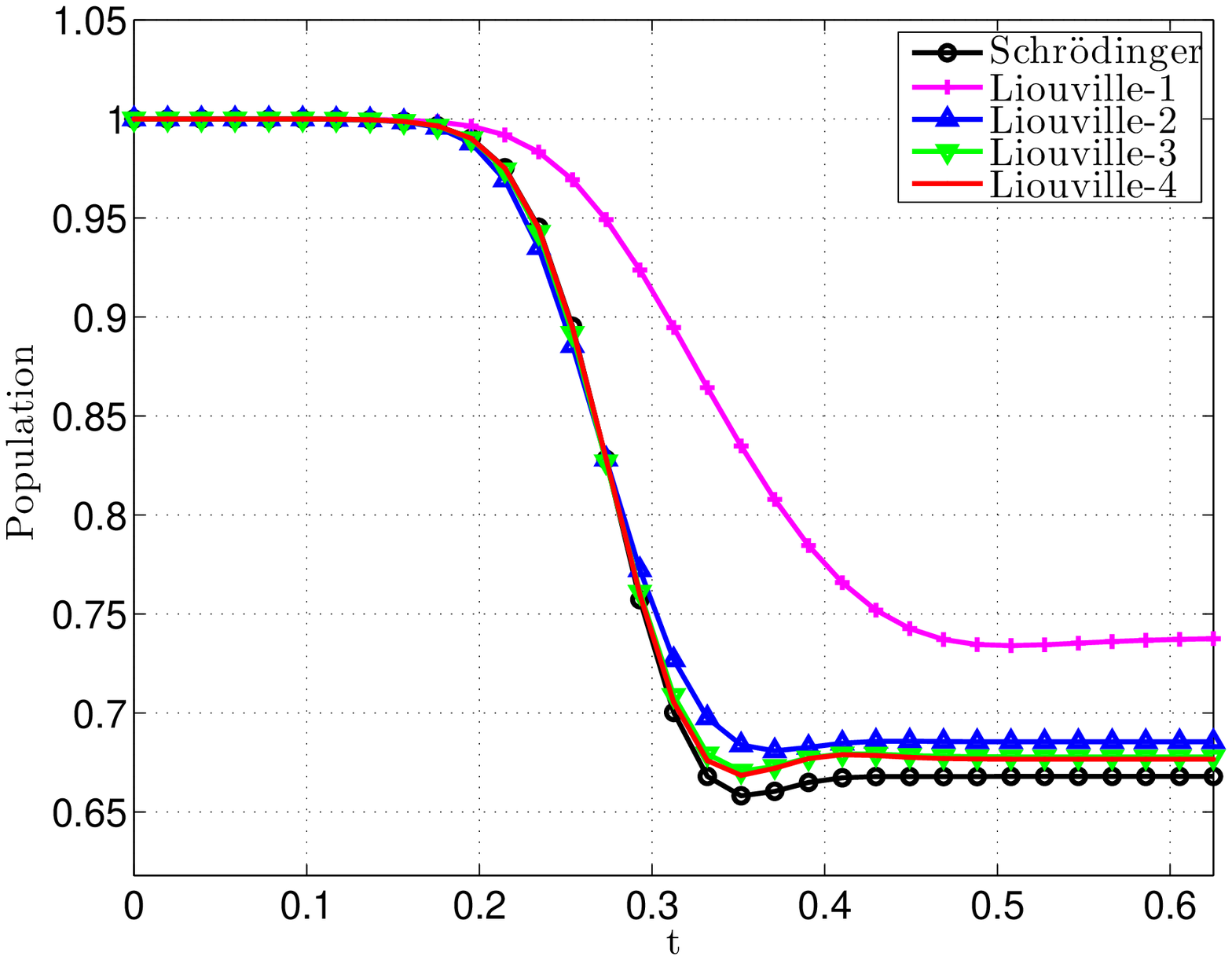}}
    \subfloat[$\varepsilon=2^{-9}$]
    {\includegraphics[width=2.5in]{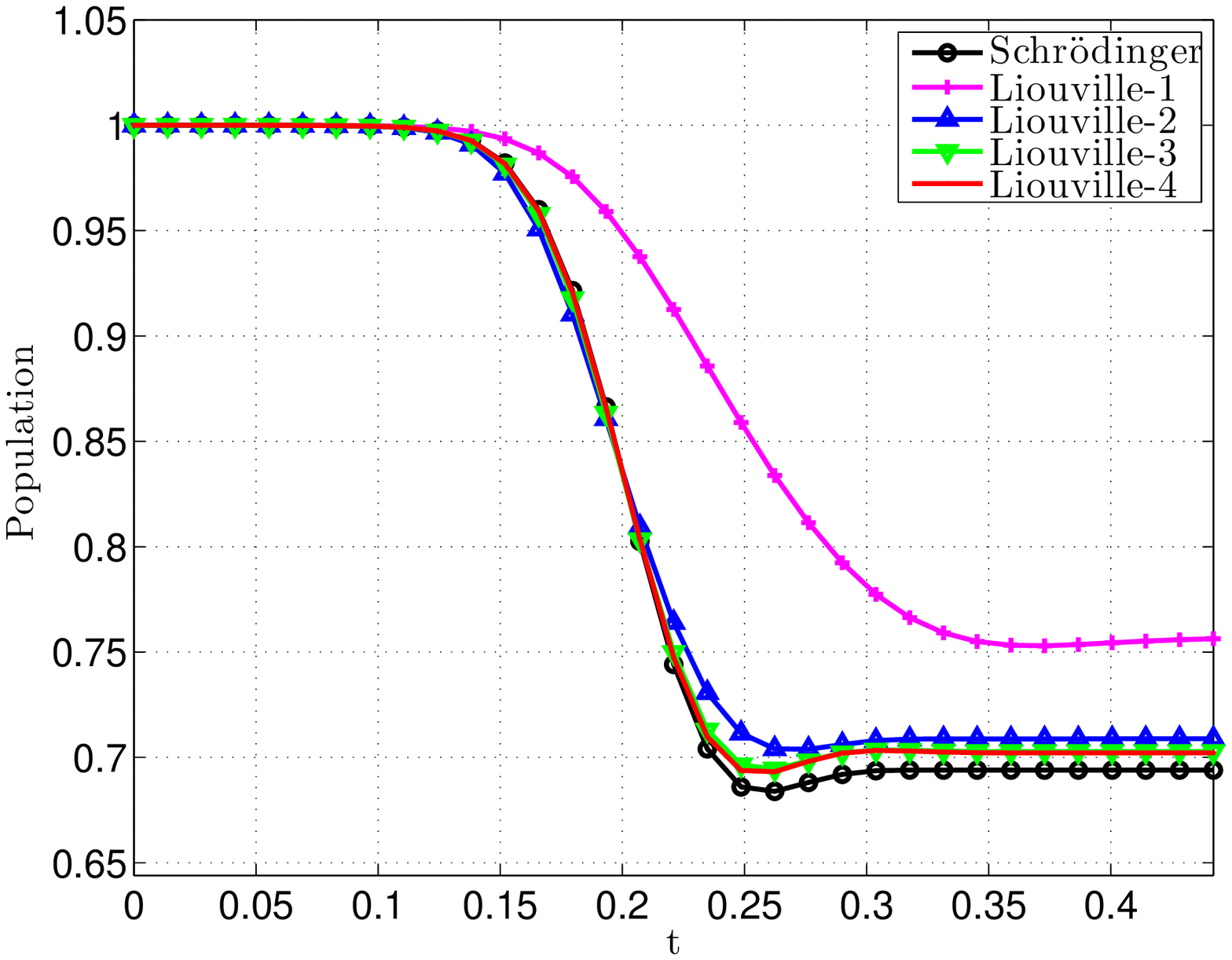}}
  \end{center}
  \caption{Two dimensional system simulation: time evolution of the population on the upper band
  $P^+_{schr/liou}$~\eqref{def:rho}. $\delta=\sqrt\vep/2$. The legend
  ``Schr\"odinger'' represents the solution of the Schr\"odinger equation,
  ``Liouville-$j$'' represents the solution of the hybrid model with $\Delta x=\Delta p = h$ in the adiabatic
  regions and $\Delta x=\Delta p = h/2$ in the non-adiabatic region, where
  $h=\sqrt{\vep}/2^{j-1}$, and $j=1,2,3,4$.}
  \label{fig:D2Pure_mass}
 \end{figure}

In this example, we deal with the problem in 2D with a pure state initial data. We set $u$ and $v$ as
\begin{equation*}\label{eqn_example_uv_2d}
 u(x)=x,\quad\text{and}\quad v(x)=\sqrt{y^2+\delta^2}.
\end{equation*}
We choose $\delta=\sqrt{\varepsilon}/2$. The minimum of energy gap is $2\delta$ and is located at the origin of the axis. The initial data for the Schr\"odinger equation are given in (\ref{eqn_initial_Sch_bands})-(\ref{eqn_initial_Sch_packet}) with $a^+=1$,
$a^-=0$, $x_0=5\sqrt{\vep}$, $y_0=0$, ${p_0}_x=-1$, and ${p_0}_y=0$. The Schr\"odinger equation is computed
using the classical time-splitting spectral method, with $\Delta x=\Delta y=\vep/8$ and
$\Delta t=5\vep^\frac{3}{2}$ and the hybrid model is computed with $\Delta x=\Delta p = h$ in the adiabatic regions
and $\Delta x=\Delta p = h/2$ in the non-adiabatic region, where
$h=O(\sqrt{\vep})$. In Figure~\ref{fig:D2Pure_density} we show the snapshots of the density contour computed by
the semiclassical model, while in Figure~\ref{fig:D2Pure_mass} we
compare the evolution of the population on the first band $P^+$ along the time given by the two systems.
In Figure~\ref{fig:D2Pure_density} we see that as time passes by, the density from the first band has some proportion jumping up to the second one. In Figure~\ref{fig:D2Pure_mass} we can clearly see that with resolved mesh the numerical solution to the semiclassical model agrees with that given by the Schr\"odinger equation.

\section{Conclusion}
In conclusion, we derived a semiclassical model for the non-adiabatic transition between different potential energy surfaces that goes beyond the classical Born-Oppenheimer approximation. By considering the complete Wigner matrix including the off-diagonal terms, our model is able to capture interesting physical phenomena such as the band-to-band transition, and the quantum correlation induced by the Berry connection. The hybrid model we proposed combines the classical adiabatic limit and the semiclassical model together to reduce the computational cost. The numerical simulations show that the hybrid model has a good agreement with the full quantum simulation.

\section*{Appendix: The Derivation of the Hamiltonian \texorpdfstring{${H'}$}{}}\label{sec_apd_system_derivation}

We give some details concerning the computation of the  ${H}$' of Eq. \eqref{eqn_H'diag}. From the definition of $\mathcal{H}$' we have
\begin{eqnarray}
  H'&=&\Tht(\bd{x})\#H(\bd{x},\bd{p})\#\Tht^\dagger(\bd{x})\nn\\
  &=&\Tht(\bd{x})\#\big(U(\bd{x})\mbb{I}+V(\bd{x})\big)\#\Tht^\dagger(\bd{x}) + \Tht(\bd{x})\#\left(\f{\bd{p}^2}{2}\mbb{I}\right)\#\Tht^\dagger(\bd{x})\nn\\
  &=&\mathrm{diag}\left\{U+E,\,U-E\right\} \,+\, \Tht(\bd{x})\#\left(\f{\bd{p}^2}{2}\mbb{I}\right)\#\Tht^\dagger(\bd{x})\nn\\
  &=&\Ld(\bd{x},\bd{p})+\ri\vep\bd{p}\cdot\nabla_{\bd{x}}\Tht(\bd{x})\Tht^\dagger(\bd{x})+\f{\vep^2}{2}\nabla_{\bd{x}}\Tht(\bd{x})\cdot\nabla_{\bd{x}}\Tht^\dagger(\bd{x}),\label{eqn_apd_H'diag}
 \end{eqnarray}
where $\Ld=U(\bd{x})+\Ld_V$ with $\Ld_V$ defined in \eqref{eqn_lambda} and we used the first order expansion of the Moyal product
\begin{eqnarray*}
  A\#B
  &=&\sum_n\f{1}{n!}\left(\f{\ri\vep}{2}\right)^nA\left(\overleftarrow\grdx\cdot\overrightarrow\grdp-\overleftarrow\grdp\cdot\overrightarrow\grdx\right)^nB\\
  &=&AB + \f{\ri\vep}{2} \left(\grdx A  \cdot \grdp B-
    \grdp A \cdot \grdx B \right)+ o (\vep).
 \end{eqnarray*}
 In particular, we used
 {
 \allowdisplaybreaks
 \beas
  \Tht(\bd{x})\#\left(\f{\bd{p}^2}{2}\mbb{I}\right)
  &=&\f{\bd{p}^2}{2}\Tht+\f{\ri\vep}{2}\bd{p}\cdot\grdx\Tht+\f{1}{2}\left(\f{\ri\vep}{2}\right)^2 \Tht\left(\overleftarrow\grdx\cdot\overrightarrow\grdp\right)^2\left(\f{\bd{p}^2}{2}\mbb{I}\right)+ o(\vep^2) \\
  &=&\f{\bd{p}^2}{2}\Tht+\f{\ri\vep}{2}\bd{p}\cdot\grdx\Tht+\f{1}{2}\left(\f{\ri\vep}{2}\right)^2 \lapx\Tht + o(\vep^2),] \eeas}
  and
  { \allowdisplaybreaks
 \beas
  \left(\f{\bd{p}^2}{2}\Tht(\bd{x})\right)\#\Tht^\dagger(\bd{x})
  &=&\f{\bd{p}^2}{2}-\f{\ri\vep}{2}\Tht\,(\bd{p}\cdot\grdx\Tht^\dagger)+\f{1}{2}\left(\f{\ri\vep}{2}\right)^2 \Tht\,\lapx\Tht^\dagger,\\
  \big(\bd{p}\cdot\grdx\Tht(\bd{x})\big)\#\Tht^\dagger(\bd{x})
  &=&\bd{p}\cdot\grdx\Tht\,\Tht^\dagger-\f{\ri\vep}{2}\grdx\Tht\cdot\grdx\Tht^\dagger,\\
  \lapx\Tht(\bd{x})\#\Tht^\dagger(\bd{x})&=&\lapx\Tht\,\Tht^\dagger, \\
  \Tht(\bd{x})\,(\bd{p}\cdot\grdx\Tht^\dagger(\bd{x}))&=&-(\bd{p}\cdot\grdx\Tht)\,\Tht^\dagger,\\
  -2\grdx\Tht(\bd{x})\cdot\grdx\Tht^\dagger(\bd{x})&=&\Tht\lapx\Tht^\dagger+\lapx\Tht\,\Tht^\dagger.
 \eeas
 }

\bibliographystyle{my-plain}

\bibliography{MyLibrary}

\end{document}